\documentstyle[12pt]{article}

\makeatletter

\def\LRA#1#2{\@tempdimb=\c@enumiv\@tempdima%
   \vcenter{\offinterlineskip\halign{##\cr%
   \hfil${\scriptstyle{#1}}$\hfil\crcr%
   \hbox to \@tempdimb{\rightarrowfill}\cr%
   \noalign{\kern-1ex}%
   \hbox to \@tempdimb{\leftarrowfill}\cr%
   \hfil${\scriptstyle{#2}}$\hfil\crcr}}}
\def\RA#1{\@tempdimb=\c@enumiv\@tempdima\vbox{\offinterlineskip%
   \halign{##\cr\hfil${\scriptstyle {#1}}$\hfil\crcr%
   \hbox to \@tempdimb{\rightarrowfill}\cr}}}
\def\LA#1{\@tempdimb=\c@enumiv\@tempdima\vbox{\offinterlineskip%
   \halign{##\cr\hfil${\scriptstyle {#1}}$\hfil\crcr%
   \hbox to \@tempdimb{\leftarrowfill}\cr}}}

\def\diag{\leavevmode\bgroup\setcounter{enumiv}{1}%
   \unitlength1em \@tempdima3em \def\\{\crcr&}\vbox\bgroup%
   \def\multicolumn##1##2{\multispan##1\setcounter{enumiv}{##1}%
   \hfil{##2}\hfil\setcounter{enumiv}{1}}
   \offinterlineskip\halign\bgroup\vrule height.8em depth.7em %
   width0pt##&&\hfil${\displaystyle{##}}$\hfil\cr&}
\def\enddiag{\crcr\egroup\egroup\egroup}
\makeatother
\if@twoside \oddsidemargin 0pt \evensidemargin 0pt \marginparwidth 85pt
\else \oddsidemargin 0pt \evensidemargin 0pt
\fi
\advance\textheight by \topskip
\font\symbolfont=msbm10
\font\teneu=eufm10
\font\egteu=eufm8
\def\dn#1{\mathchoice{\hbox{\teneu #1}}{\hbox{\teneu #1}}%
   {\hbox{\egteu #1}}{\hbox{\egteu #1}}}

\def\ccc{\hbox{\symbolfont C}}
\def\fff{\hbox{\symbolfont F}}

\def\ggg{\hbox{\symbolfont G}}
\def\hhh{\hbox{\symbolfont H}}
\def\mmm{\hbox{\symbolfont M}}

\def\zzz{\hbox{\symbolfont Z}}
\def\ppp{\hbox{\symbolfont P}}
\def\qqq{\hbox{\symbolfont Q}}

\def\depth{\mathop d\nolimits}

\def\aa{\mathop{\cal A}\nolimits}
\def\bb{\mathop{\cal B}\nolimits}
\def\cc{\mathop{\cal C}\nolimits}
\def\dd{\mathop{\cal D}\nolimits}

\def\ii{\mathop{\cal I}\nolimits}
\def\pp{\mathop{\cal P}\nolimits}

\def\tt{\mathop{\cal T}\nolimits}
\def\xx{\mathop{\cal X}\nolimits}
\def\yy{\mathop{\cal Y}\nolimits}

\def\pd{\mathop{\rm pd}\nolimits}
\def\height{\mathop{\rm ht}\nolimits}
\def\id{\mathop{\rm id}\nolimits}
\def\gl{\mathop{\rm gl.dim}\nolimits}
\def\rdim{\mathop{\rm rep.dim}\nolimits}

\def\height{\mathop{\rm ht}\nolimits}
\def\Spec{\mathop{\rm Spec}\nolimits}

\def\spep#1{\mathop{{}^{\bullet}\strut\kern-.1em{#1}}\nolimits}

\def\hom{\mathop{\rm Hom}\nolimits}
\def\endm{\mathop{\rm End}\nolimits}
\def\aut{\mathop{\rm Aut}\nolimits}
\def\inn{\mathop{\rm Inn}\nolimits}
\def\out{\mathop{\rm Out}\nolimits}
\def\add{\mathop{\rm add}\nolimits}
\def\Mod{\mathop{\rm Mod}\nolimits}

\def\ind{\mathop{\rm ind}\nolimits}
\def\irr{\mathop{\rm irr}\nolimits}
\def\cm{\mathop{\rm CM}\nolimits}

\def\Cok{\mathop{\rm Cok}\nolimits}
\def\Ker{\mathop{\rm Ker}\nolimits}
\def\cen{\mathop{\rm Cen}\nolimits}

\def\der#1{\mathop{{\rm R}^{#1}\alpha}}

\def\mod{\mathop{\rm mod}\nolimits}

\def\cen{\mathop{\rm Cen}\nolimits}
\def\add{\mathop{\rm add}\nolimits}

\def\cc{\mathop{\cal C}\nolimits}

\def\Mod{\mathop{\rm Mod}\nolimits}
\def\Cok{\mathop{\rm Cok}\nolimits}
\def\Ker{\mathop{\rm Ker}\nolimits}

\def\depth{\mathop{\rm depth}\nolimits}
\def\pd{\mathop{\rm pd}\nolimits}
\def\fd{\mathop{\rm fd}\nolimits}
\def\id{\mathop{\rm id}\nolimits}
\def\gl{\mathop{\rm gl.dim}\nolimits}
\def\domdim{\mathop{\rm dom.dim}\nolimits}

\def\height{\mathop{\rm ht}\nolimits}
\def\Spec{\mathop{\rm Spec}\nolimits}

\def\Gl{\mathop{\rm GL}\nolimits}

\def\ma{\mathop{\rm M}\nolimits}
\def\ext{\mathop{\rm Ext}\nolimits}

\def\grade{\mathop{{\rm grade}}}

\def\resdim#1#2{\mathop{{#1}\mbox{-{\rm dim}}\strut\kern2pt {#2}}\nolimits}

\def\timesl{\begin{picture}(16,8)\put(3,-1){\line(1,1){10}}\put(13,9){\line(0,-1){10}}\put(13,-1){\line(-1,1){10}}\end{picture}}
\def\XZA{0.1}
\def\XZB{0.2}
\def\XA{1}
\def\XAA{1.1}
\def\XAB{1.2}
\def\XB{2}
\def\XBA{2.1}
\def\XBB{2.2}
\def\XBBA{2.2.1}
\def\XBBB{2.2.2}
\def\XBBC{2.2.3}
\def\XBC{2.3}
\def\XBCA{2.3.1}
\def\XBCB{2.3.2}

\def\XBD{2.4}
\def\XBDA{2.4.1}
\def\XBE{2.5}
\def\XBEA{2.5.1}
\def\XBEB{2.5.2}
\def\XBEC{2.5.3}
\def\XBF{2.6}
\def\XBFA{2.6.1}
\def\XBFB{2.6.2}
\def\XC{3}
\def\XCA{3.1}
\def\XCAA{3.1.1}
\def\XCAB{3.1.2}
\def\XCB{3.2}
\def\XCBA{3.2.1}
\def\XCBB{3.2.2}
\def\XCBC{3.2.3}
\def\XCC{3.3}
\def\XCCA{3.3.1}
\def\XCD{3.4}
\def\XCDA{3.4.1}
\def\XCDB{3.4.2}
\def\XCDC{3.4.3}
\def\XCE{3.5}
\def\XCEA{3.5.1}
\def\XCEB{3.5.2}
\def\XCEC{3.5.3}
\def\XD{4}
\def\XDA{4.1}
\def\XDAA{4.1.1}
\def\XDB{4.2}
\def\XDBA{4.2.1}
\def\XDBB{4.2.2}
\def\XDBC{4.2.3}
\def\XDBD{4.2.4}
\def\XDC{4.3}
\def\XDCA{4.3.1}
\def\XDD{4.4}
\def\XDDA{4.4.1}
\def\XDDB{4.4.2}
\def\XDDC{4.4.3}
\def\XDDD{4.4.4}
\def\XDE{4.5}
\def\XDF{4.6}
\def\XDG{4.7}
\def\XDGA{4.7.1}

\def\XDH{4.8}
\def\XE{5}
\def\XEA{5.1}
\def\XEB{5.2}
\def\XEBA{5.2.1}
\def\XEBB{5.2.2}
\def\XEC{5.3}
\def\XECA{5.3.1}
\def\XECB{5.3.2}
\def\XECC{5.3.3}
\def\XECD{5.3.4}
\def\XED{5.4}
\def\XEDA{5.4.1}
\def\XEDB{5.4.2}
\def\XEDC{5.4.3}
\def\XEDD{5.4.4}
\def\XEE{5.5}
\def\XEEA{5.5.1}
\def\XEEB{5.5.2}
\def\XF{6}
\def\XFA{6.1}
\def\XFAA{6.1.1}
\def\XFAB{6.1.2}
\def\XFAC{6.1.3}
\def\XFAD{6.1.4}
\def\XFB{6.2}
\def\XFBA{6.2.1}
\def\XFBB{6.2.2}
\def\XFBC{6.2.3}
\def\XFBD{6.2.4}
\def\XFC{6.3}
\def\XFCA{6.3.1}
\begin{document}
\begin{center}
\vspace*{1cm}{\large Auslander correspondence}
\vskip1em Osamu Iyama
\end{center}

{\footnotesize
{\sc Abstract. }We study Auslander correspondence from the viewpoint of higher dimensional Auslander-Reiten theory on maximal orthogonal subcategories. We give homological characterizations of Auslander algebras, especially an answer to a question of M. Artin. They are also closely related to Auslander's representation dimension of artin algebras and Van den Bergh's non-commutative crepant resolutions of Gorenstein singularities.}

\vskip.5em
Let us recall M. Auslander's classical theorem [A1] below, which introduced a completely new insight to the representation theory of algebras.

\vskip.5em{\bf\XZA\ Theorem }(Auslander correspondence) {\it
There exists a bijection between the set of Morita-equivalence classes of representation-finite finite-dimensional algebras $\Lambda$ and that of finite-dimensional algebras $\Gamma$ with $\gl\Gamma\le2$ and $\domdim\Gamma\ge2$ (\XCC). It is given by $\Lambda\mapsto\Gamma:=\endm_\Lambda(M)$ for an additive generator $M$ of $\mod\Lambda$.}

\vskip.5em
In this really surprising theorem, the representation theory of $\Lambda$ is encoded in the structure of the homologically nice algebra $\Gamma$ called an {\it Auslander algebra}. Since the category $\mod\Gamma$ is equivalent to the functor category on $\mod\Lambda$, Auslander correspondence gave us a prototype of the use of functor categories in representation theory. In this sense, Auslander correspondence was a starting point of later Auslander-Reiten theory [ARS] historically. Theoretically, Auslander correspondence gives a direct connection between two completely different concepts, i.e. a representation theoretic property `representation-finiteness' and a homological property `$\gl\Gamma\le2$ and $\domdim\Gamma\ge2$'. It is a quite interesting project to find correspondence between representation theoretic properties and homological properties. Algebras $\Gamma$ with $\gl\Gamma\le2$ was studied in [I4,5] from this viewpoint.

The aim of this paper is to give a higher dimensional version of Auslander correspondence. Recently, it was pointed out in [I7] that Auslander-Reiten theory is `2-dimensional-like' theory, and the concept of maximal $(n-1)$-orthogonal subcategories was introduced as a domain of `$(n+1)$-dimensional' Auslander-Reiten theory. Thus it would be natural to study `$(n+1)$-dimensional' Auslander correspondence from the viewpoint of maximal $(n-1)$-orthogonal subcategories. One of our main results is the theorem below, which is a special case of \XDBB. We call an additive category {\it finite} if it has an additive generator.

\vskip.5em{\bf\XZB\ Theorem }($(n+1)$-dimensional Auslander correspondence) {\it
For any $n\ge1$, there exists a bijection between the set of equivalence classes of finite maximal $(n-1)$-orthogonal subcategories $\cc$ of $\mod\Lambda$ for finite-dimensional algebras $\Lambda$, and the set of Morita-equivalence classes of finite-dimensional algebras $\Gamma$ with $\gl\Gamma\le n+1$ and $\domdim\Gamma\ge n+1$. It is given by $\cc\mapsto\Gamma:=\endm_\Lambda(M)$ for an additive generator $M$ of $\cc$.}

\vskip.5em
Putting $n=1$ in this theorem, we obtain the theorem \XZA\ because $\mod\Lambda$ has a unique maximal $0$-orthogonal subcategory $\mod\Lambda$.

We study not only a higher global-dimensional version of \XZA\ but also its higher Krull-dimensional version. Auslander-Reiten theory plays a crucial role also in the study of the category $\cm\Lambda$ of Cohen-Macaulay modules over Cohen-Macaulay rings and orders $\Lambda$ (\XCA) of Krull-dimension $d$ [A2,3,4][AR2][Y]. A version of \XZA\ for the case $d=1$ and $2$ was given by Auslander-Roggenkamp [ARo] and Auslander [Ar][RV] respectively. But it seems that any version of \XZA\ for the case $d>2$ is unknown. Especially, M. Artin raised a question in [Ar] to characterize homologically the endomorphism rings $\endm_\Lambda(M)$ of an additive generator $M$ of $\cm\Lambda$ for representation-finite orders $\Lambda$ with $d>2$. In \XDBC, we give an answer to this question.

Since the category $\cm\Lambda$ over an $R$-order $\Lambda$ is the orthogonal category ${}^\perp T$ for $T:=\hom_R(\Lambda,R)$, we study the orthogonal category $\bb:={}^\perp T$ for cotilting $\Lambda$-modules $T$ with $\id{}_\Lambda T=m$ (\XCB) in general. The category $\bb$ seems to be still `2-dimensional-like' even if $m>2$ from the viewpoint of [I7], and we study `$(n+1)$-dimensional' Auslander-Reiten theory on maximal $(n-1)$-orthogonal subcategories $\cc$ of $\bb$ in \S\XB. We call the endomorphism ring $\endm_\Lambda(M)$ of an additive generator $M$ of $\cc$ an {\it Auslander algebra of type $(d,m,n)$}, and give homological characterizations in \S\XD. We see in \XDBB\ that Auslander correspondence of type $(d,m,n)$ can be stated in terms of the {\it $(m+1,n+1)$-condition} (\XCC) which was introduced in [I2,4] as a bridge between Auslander's $n$-Gorenstein condition [FGR][B][AR4][C] and the dominant dimension [T][H]. Since `higher dimensional' Auslander-Reiten theory for the case $d=m=n+1$ is quite peculiar [I7], Auslander algebras of type $(d,d,d-1)$ have a very nice homological characterizations in \XDG, especially (3) is closely related to Artin-Schelter regular ring of dimension $d$. We observe in \XDGA\ that the $(n+1,n+1)$-condition means the existence of $n$-almost split sequence homologically.

Recently, in representation theory and non-commutative algebraic geometry, it seems that the study of `nice' subcategories becomes more and more important besides our maximal $(n-1)$-orthogonal subcategories. Especially, Van den Bergh introduced the concept of non-commutative crepant resolutions [V1,2] to study the Bondal-Orlov conjecture [BO] on derived categories of resolutions of a Gorenstein singularity. We see in \XEBA\ that non-commutative crepant resolutions are almost one same concept as our maximal $(d-1)$-orthogonal subcategories, and in \XECC\ that all maximal $1$-orthogonal subcategories are derived equivalent, which supports Van den Bergh's generalization [V2] of the Bondal-Orlov conjecture. As was pinted out by Leuschke [L], we see in \XED\ that the concept of non-commutative crepant resolutions is also closely related to the concept of Auslander's representation dimension [A1], which measures how far an algebra is from being representation-finite. A lot of recent results on the representation dimension (see references in \XEDD(1)) show that it is a really interesting and useful concept. Although the representation dimension is always finite for $d\le1$ [I1,3,6], we see in \XEDC\ that this is not the case for $d\ge2$. We give in \XEE\ a boundedness conjecture for $1$-orthogonal subcategories, and prove it for algebras with the representation dimension at most three.

In \S\XF, we give three remarkable examples. In \XFA, we observe a higher dimensional version of Auslander's theorem on McKay correspondence [A4]. In \XFB, we see that the work of Geiss-Leclerc-Schr\"oer [GLS1,2] on rigid modules on preprojective algebras is closely related to our study. In \XFC, we see that the work of Buan-Marsh-Reineke-Reiten-Todorov [BMRRT] on cluster categories is also closely related to our study.

\vskip-.5em
\begin{center}
{\sc Acknowledgements}
\end{center}

\vskip-.5em
The author would like to thank Idun Reiten, Jan Schr\"oer and Yuji Yoshino for stimulating discussion. He thanks Idun Reiten for her kind advice in correcting my English. He also would like to thank Michel Van den Bergh who suggested \XEBA\ for a certain case.

\vskip.5em{\bf\XA\ Preliminaries on functor categories }

\vskip.5em{\bf\XAA\ }
Let $\aa$ be an additive category and $\cc$ a subcategory of $\aa$.

(1) We denote by $\aa(X,Y)$ the set of morphisms from $X$ to $Y$, and by $fg\in\aa(X,Z)$ the composition of $f\in\aa(X,Y)$ and $g\in\aa(Y,Z)$. We denote by $J_{\aa}$ the Jacobson radical of $\aa$, and by $\ind\aa$ the set of isoclasses of indecomposable objects in $\aa$. An {\it $\aa$-module} is a contravariant additive functor from $\aa$ to the category of abelian groups. We denote by $\Mod\aa$ the abelian category of $\aa$-modules. We call an $\aa$-module $F$ {\it finitely presented} if there exists an exact sequence $\aa(\ ,X)\stackrel{\cdot f}{\to}\aa(\ ,Y)\to F\to0$. We denote by $\mod\aa$ the category of finitely presented $\aa$-modules [AR1]. We call $f\in J_{\aa}(X,Y)$ a {\it sink map} if $\aa(\ ,X)\stackrel{\cdot f}{\to}J_{\aa}(\ ,Y)\to0$ is exact, and a {\it source map} if $\aa(Y,\ )\stackrel{f\cdot}{\to} J_{\aa}(X,\ )\to0$ is exact.

We call $\cc$ {\it contravariantly} (resp. {\it covariantly}) {\it finite} in $\aa$ if $\aa(\ ,X)|_{\cc}$ (resp. $\aa(X,\ )|_{\cc}$) is a finitely generated $\cc$-module for any $X\in\aa$ [AS]. We call $\cc$ {\it functorially finite} if it is contravariantly and covariantly finite. We call a complex $\cdots\stackrel{f_2}{\to}C_1\stackrel{f_1}{\to}C_0\stackrel{f_0}{\to}X$ a {\it right $\cc$-resolution} of $X\in\aa$ if $C_i\in\cc$ and $\cdots\stackrel{\cdot f_2}{\to}\aa(\ ,C_1)\stackrel{\cdot f_1}{\to}\aa(\ ,C_0)\stackrel{\cdot f_0}{\to}\aa(\ ,X)\to0$ is exact on $\cc$. We write $\resdim{\cc}{X}\le n$ if $X$ has a right $\cc$-resolution with $C_{n+1}=0$. Put $\resdim{\cc}{\aa}:=\sup_{X\in\aa}\resdim{\cc}{X}$. Define a {\it left $\cc$-resolution}, $\resdim{\cc{}^{op}}{X}$ and $\resdim{\cc{}^{op}}{\aa^{op}}$ dually. We denote by $[\cc]$ the ideal of $\aa$ consisting of morphisms which factor through $\cc$. 

(2) Let $R$ be a commutative local ring and $D:\Mod R\to\Mod R$ the Matlis dual. Assume that $\aa$ is an $R$-category such that $\aa(X,Y)$ is an $R$-module of finite length for any $X,Y\in\aa$. For any $F\in\Mod\aa$ and $X\in\aa$, $F(X)$ has an $R$-module structure naturally. Thus we have a functor $D:\Mod\aa\leftrightarrow\Mod\aa^{op}$ by composing with $D$. We call $\aa$ a {\it dualizing $R$-variety} if $D$ induces a duality $\mod\aa\leftrightarrow\mod\aa^{op}$ [AR1]. If $\aa$ is a dualizing $R$-variety, it is easily checked that $\mod\aa$ (resp. $\mod\aa^{op}$) is an abelian subcategory of $\Mod\aa$ (resp. $\Mod\aa^{op}$) which is closed under kernels, cokernels and extensions [A1;2.1]. In particular, $\aa$ has pseudo-kernels and pseudo-cokernels.

\vskip.5em{\bf\XAB\ }The following version of a theorem of Auslander-Smalo [AS;2.3] gives a relationship between dualizing $R$-varieties and functorially finite subcategories.

\vskip.5em{\bf Proposition }{\it
Let $\aa$ be a dualizing $R$-variety. Then any functorially finite subcategory $\cc$ of $\aa$ is a dualizing $R$-variety.}

\vskip.5em{\sc Proof }
(i) We will show that $F|_{\cc}\in\mod\cc$ holds for any $F\in\mod\aa$.

Since $\mod\cc$ is closed under cokernels in general, we only have to show that $\aa(\ ,X)|_{\cc}\in\mod\cc$ holds for any $X\in\aa$. Let $f\in\aa(C_0,X)$ be a right $\cc$-resolution, $g\in\aa(X_1,C_0)$ a pseudo-kernel of $f$, and $h\in\aa(C_1,X_1)$ a right $\cc$-resolution. Then $\aa(\ ,C_1)\stackrel{\cdot hg}{\to}\aa(\ ,C_0)\stackrel{\cdot f}{\to}\aa(\ ,X)\to0$ is exact on $\cc$.

(ii) For any $F\in\mod\cc$, take an exact sequence $\cc(\ ,Y)\stackrel{\cdot f}{\to}\cc(\ ,X)\to F\to0$. Define $F^\prime\in\mod\aa$ by an exact sequence $\aa(\ ,Y)\stackrel{\cdot f}{\to}\aa(\ ,X)\to F^\prime\to0$. Since $\aa$ is a dualizing $R$-variety, $DF^\prime\in\mod\aa^{op}$ holds. Thus $DF=(DF^\prime)|_{\cc}\in\mod\cc^{op}$ holds by (i). A dual argument shows that $DG\in\mod\cc$ holds for any $G\in\mod\cc{}^{op}$.\rule{5pt}{10pt}

\vskip.5em{\bf\XB\ Higher dimensional Auslander-Reiten theory for dualizing $R$-varieties}

In Auslander-Reiten theory, there are two approaches to showing the existence theorem of almost split sequences. One is based on an explicit calculation of extension groups [ARS], and higher dimensional Auslander-Reiten theory in [I7] was developped in this direction. Another is more general and suggestive but less concrete, and based on the concept of dualizing $R$-varieties [AR1][AS]. In this section, we will study higher dimensional Auslander-Reiten theory in the latter direction. This will enable us to treat the orthogonal category ${}^\perp T$ for a cotilting $\Lambda$-module $T$ in \S\XC.

\vskip.5em{\bf\XBA\ }
Throughout this section, assume that $R$ is a commutative local ring and $\aa$ is an abelian $R$-category with enough projectives. For $X,Y\in\aa$, we write $X\perp_n Y$ if $\ext^i_{\aa}(X,Y)=0$ holds for any $i$ ($0<i\le n$). Put $\cc^{\perp_n}:=\{X\in\aa\ |\ \cc\perp_n X\}$ and ${}^{\perp_n}\cc:=\{X\in\aa\ |\ X\perp_n\cc\}$. Put $\perp:=\perp_\infty$. Let $\pp=\pp(\aa):={}^\perp\aa$ be the category of projective objects in $\aa$. Let $\underline{\aa}:=\aa/[\pp]$ be the {\it stable category} (\XAA), and $\Omega:\underline{\aa}\to\underline{\aa}$ the {\it syzygy functor}. One can easily check the facts below for any $X\in{}^{\perp_n}\pp$ and $Y\in\pp$.

(1) $\Omega^n:\underline{\aa}(X,Y)\to\underline{\aa}(\Omega^nX,\Omega^nY)$ is bijective.

(2) We have a functorial isomorphism $\underline{\aa}(\Omega^nX,Y)=\ext^n_{\aa}(X,Y)$.

\vskip.5em{\bf\XBB\ }
In the rest of this section, we assume that $\bb$ is a {\it resolving} subcategory of $\aa$, i.e. $\pp\subseteq\bb$ and $\bb$ is closed under extensions and kernels of surjections [AR3]. Thus $\Omega$ induces the syzygy functor $\Omega:\underline{\bb}\to\underline{\bb}$ for $\underline{\bb}:=\bb/[\pp]$. Let $\ii=\ii(\bb):=\bb^\perp\cap\bb$ be the category of injective objects in $\bb$. Moreover, we assume that $\bb$ is {\it enough injectives}, i.e. for any $X\in\bb$, there exists an exact sequence $0\to X\to I\to Y\to0$ with $Y\in\bb$ and $I\in\ii$. Let $\overline{\bb}:=\bb/[\ii]$ be the {\it costable category}, and $\Omega^-:\overline{\bb}\to\overline{\bb}$ the {\it cosyzygy functor}. For a subcategory $\cc$ of $\bb$, we denote by $\underline{\cc}$ (resp. $\overline{\cc}$) the corresponding subcategory of $\underline{\bb}$ (resp. $\overline{\bb}$). It is not difficult to check the proposition below (cf. [AR1]).

\vskip.5em{\bf\XBBA\ Proposition }{\it
Let $0\to X_2\stackrel{g}{\to}X_1\stackrel{f}{\to}X_0\to0$ be an exact sequence in $\aa$ with $X_i\in\bb$. Then we have the two long exact sequences below.
\begin{eqnarray*}
&\cdots\stackrel{}{\to}\underline{\bb}(\ ,\Omega X_2)\stackrel{\cdot\Omega g}{\to}\underline{\bb}(\ ,\Omega X_1)\stackrel{\cdot\Omega f}{\to}\underline{\bb}(\ ,\Omega X_0)\stackrel{}{\to}\underline{\bb}(\ ,X_2)\stackrel{\cdot g}{\to}\underline{\bb}(\ ,X_1)\stackrel{\cdot f}{\to}\underline{\bb}(\ ,X_0)&\\
&\stackrel{}{\to}\ext^1_{\aa}(\ ,X_2)\stackrel{\cdot g}{\to}\ext^1_{\aa}(\ ,X_1)\stackrel{\cdot f}{\to}\ext^1_{\aa}(\ ,X_0)\stackrel{}{\to}\ext^2_{\aa}(\ ,X_2)\stackrel{\cdot g}{\to}\ext^2_{\aa}(\ ,X_1)\stackrel{\cdot f}{\to}\cdots&
\end{eqnarray*}
\vskip-1.3em\begin{eqnarray*}
&\cdots\stackrel{}{\to}\overline{\bb}(\Omega^-X_0,\ )\stackrel{\Omega^-f\cdot}{\to}\overline{\bb}(\Omega^-X_1,\ )\stackrel{\Omega^-g\cdot}{\to}\overline{\bb}(\Omega^-X_2,\ )\stackrel{}{\to}\overline{\bb}(X_0,\ )\stackrel{f\cdot}{\to}\overline{\bb}(X_1,\ )\stackrel{g\cdot}{\to}\overline{\bb}(X_2,\ )&\\
&\stackrel{}{\to}\ext^1_{\aa}(X_0,\ )\stackrel{f\cdot}{\to}\ext^1_{\aa}(X_1,\ )\stackrel{g\cdot}{\to}\ext^1_{\aa}(X_2,\ )\stackrel{}{\to}\ext^2_{\aa}(X_0,\ )\stackrel{f\cdot}{\to}\ext^2_{\aa}(X_1,\ )\stackrel{g\cdot}{\to}\cdots.&
\end{eqnarray*}}

\vskip-1em{\bf\XBBB\ }The following fundamental theorem is a version of [AR1].

\vskip.5em{\bf Theorem }{\it
(1) $\mod\underline{\bb}$ is enough injectives, and $X\mapsto\ext^1_{\aa}(\ ,X)$ gives an equivalence from $\overline{\bb}$ to the category $\ii(\mod\underline{\bb})$ of finitely presented injective $\underline{\bb}$-modules.

(2) $\mod\overline{\bb}^{op}$ is enough injectives, and $X\mapsto\ext^1_{\aa}(X,\ )$ gives an equivalence from $\underline{\bb}^{op}$ to the category $\ii(\mod\overline{\bb}^{op})$ of finitely presented injective $\overline{\bb}^{op}$-modules.}

\vskip.5em{\sc Proof }We only prove (1) since (2) is proved dually. Fix $X,Y\in\bb$. Let $0\to X\stackrel{g}{\to}I\stackrel{f}{\to}\Omega^-X\to0$ be an injective resolution. Then ${\bf A}:\underline{\bb}(\ ,I)\stackrel{\cdot f}{\to}\underline{\bb}(\ ,\Omega^-X)\to\ext^1_{\aa}(\ ,X)\to0$ is exact by \XBBA. Thus $\ext^1_{\aa}(\ ,X)\in\mod\underline{\bb}$. We have an exact sequence $\hom({\bf A},\ext^1_{\aa}(\ ,Y)):0\to\hom(\ext^1_{\aa}(\ ,X),\ext^1_{\aa}(\ ,Y))\to\ext^1_{\aa}(\Omega^-X,Y)\stackrel{f\cdot}{\to}\ext^1_{\aa}(I,Y)$ by Yoneda's lemma. Since $0\to\overline{\bb}(X,Y)\to\ext^1_{\aa}(\Omega^-X,Y)\stackrel{f\cdot}{\to}\ext^1_{\aa}(I,Y)$ is exact by \XBBA, we have a bijection $\overline{\bb}(X,Y)\to\hom(\ext^1_{\aa}(\ ,X),\ext^1_{\aa}(\ ,Y))$. Thus the functor $\overline{\bb}\to\mod\underline{\bb}$ given by $X\mapsto\ext^1_{\aa}(\ ,X)$ is full and faithful.

For any $F\in\mod\underline{\bb}$, take an exact sequence $\underline{\bb}(\ ,Y_1)\stackrel{\cdot f}{\to}\underline{\bb}(\ ,Y_0)\to F\to0$. Then $f$ is an epimorphism since $F(\pp)=0$. Let $0\to Y_2\stackrel{g}{\to}Y_1\stackrel{f}{\to}Y_0\to0$ be an exact sequence in $\aa$. Then $Y_2\in\bb$ (Be careful in the proof of (2)). By \XBBA, ${\bf P}:\underline{\bb}(\ ,Y_2)\stackrel{\cdot g}{\to}\underline{\bb}(\ ,Y_1)\stackrel{\cdot f}{\to}\underline{\bb}(\ ,Y_0)\to F\to0$ gives a projective resolution of $F$. We have an exact sequence $\hom({\bf P},\ext^1_{\aa}(\ ,X)):\ext^1_{\aa}(Y_0,X)\stackrel{f\cdot}{\to}\ext^1_{\aa}(Y_1,X)\stackrel{g\cdot}{\to}\ext^1_{\aa}(Y_2,X)$ by Yoneda's lemma and \XBBA. Thus $\ext^1(F,\ext^1_{\aa}(\ ,X))=0$ holds, and $\ext^1_{\aa}(\ ,X)$ is injective. Since we have an exact sequence $0\to F\to\ext^1_{\aa}(\ ,Y_2)$ by \XBBA, $\mod\underline{\bb}$ is enough injectives.\rule{5pt}{10pt}

\vskip.5em{\bf\XBBC\ }
In the rest of this section, we assume that the conditions in the following version of a theorem of Auslander-Reiten [AR5;2.2] are satisfied.

\vskip.5em{\bf Proposition }{\it
Assume that $\ext^1_{\aa}(X,Y)$ is an $R$-module of finite length for any $X,Y\in\bb$. Then the following conditions are equivalent.

(1) $\underline{\bb}$ is a dualizing $R$-variety.

(2) $\overline{\bb}$ is a dualizing $R$-variety.

(3) There exists an equivalence $\tau:\underline{\bb}\to\overline{\bb}$ with a quasi-inverse $\tau^-$ and a functorial isomorphism $\overline{\bb}(Y,\tau X)\simeq D\ext^1_{\aa}(X,Y)\simeq \underline{\bb}(\tau^-Y,X)$ for any $X,Y\in\bb$.}

\vskip.5em{\sc Proof }
(2)$\Rightarrow$(3) By Yoneda's lemma, $X\mapsto\overline{\bb}(\ ,X)$ gives an equivalence $\fff:\overline{\bb}\to\pp(\mod\overline{\bb})$. By (2) and \XBBB, $X\mapsto D\ext^1_{\aa}(X,\ )$ gives an equivalence $\ggg:\underline{\bb}\to\pp(\mod\overline{\bb})$. Let $\fff^-$ be a quasi-inverse of $\fff$, $\tau:=\fff^-\circ\ggg$ and $\tau^-$ a quasi-inverse of $\tau$. Then the assertion follows. A dual argument shows (1)$\Rightarrow$(3).

(3)$\Rightarrow$(1)$\land$(2) Fix $F\in\mod\underline{\bb}$. By \XBBA, we can take an exact sequence $0\to F\to\ext^1_{\aa}(\ ,X_2)\stackrel{\cdot g}{\to}\ext^1_{\aa}(\ ,X_1)$, which is induced by an exact sequence $0\to X_2\stackrel{g}{\to}X_1\stackrel{f}{\to}X_0\to0$ in $\bb$. Applying $D$, we have an exact sequence $\underline{\bb}(\tau^-X_1,\ )\stackrel{\tau^-g\cdot}{\longrightarrow}\underline{\bb}(\tau^-X_2,\ )\to DF\to0$ by (3). Thus $DF\in\mod\underline{\bb}^{op}$ holds. Dually, $DG\in\mod\overline{\bb}$ holds for any $G\in\mod\overline{\bb}^{op}$. Since we have equivalences $\tau:\mod\overline{\bb}\to\mod\underline{\bb}$ and $\tau:\mod\overline{\bb}^{op}\to\mod\underline{\bb}^{op}$ which commute with $D$, we obtain (1) and (2).\rule{5pt}{10pt}

\vskip.5em{\bf\XBC\ }For $n\ge1$, we define functors $\tau_n$ and $\tau_n^-$ by \[\tau_n:=\tau\circ\Omega^{n-1}:\underline{\bb}\to\overline{\bb}\ \ \ \ \ {\rm and }\ \ \ \ \ \tau^-_n:=\tau^-\circ\Omega^{-(n-1)}:\overline{\bb}\to\underline{\bb},\]
where $\Omega:\underline{\bb}\to\underline{\bb}$ is the syzygy functor and $\Omega^{-}:\overline{\bb}\to\overline{\bb}$ is the cosyzygy functor. Put
\[\xx_n:={}^{\perp_{n}}\pp\cap\bb\ \ \ \ \ {\rm and }\ \ \ \ \ \yy_n:=\ii^{\perp_{n}}\cap\bb.\]
Let us give a version of [I7;1.4,1.5] for our situation.

\vskip.5em{\bf\XBCA\ Theorem }{\it
(1) There exist functorial isomorphisms $\overline{\bb}(Y,\tau_nX)\simeq D\ext^n_{\aa}(X,Y)\simeq\underline{\bb}(\tau_n^-Y,X)$ for any $X,Y\in\bb$. Thus $\tau_n:\underline{\bb}\to\overline{\bb}$ is a right adjoint of $\tau_n^-:\overline{\bb}\to\underline{\bb}$.

(2) For any $i$ ($0<i<n$), there exist functorial isomorphisms below for any $X\in\xx_{n-1}$, $Y\in\yy_{n-1}$ and $Z\in\bb$.
{\small\begin{eqnarray*}
&D\ext^n_{\aa}(X,Z)\simeq\overline{\bb}(Z,\tau_nX),\ \ D\ext^{n-i}_{\aa}(X,Z)\simeq\ext^{i}_{\aa}(Z,\tau_nX),\ \ D\underline{\bb}(X,Z)\simeq\ext^n_{\aa}(Z,\tau_nX)&\\
&D\ext^n_{\aa}(Z,Y)\simeq\underline{\bb}(\tau_n^-Y,Z),\ \ D\ext^{n-i}_{\aa}(Z,Y)\simeq\ext^{i}_{\aa}(\tau_n^-Y,Z),\ \ D\overline{\bb}(Z,Y)\simeq\ext^n_{\aa}(\tau_n^-Y,Z)&
\end{eqnarray*}}}

\vskip-1em{\sc Proof }
(1) We have functorial isomorphisms $\overline{\bb}(Y,\tau_nX)\stackrel{\XBBC(3)}{\simeq}D\ext^1_{\aa}(\Omega^{n-1}X,Y)\stackrel{}{\simeq}D\ext^n_{\aa}(X,Y)$. The isomorphism for $\tau^-_n$ is given dually.

(2) The left isomorphisms are given in (1). For $i>0$, we have functorial isomorphisms $\ext^i_{\aa}(Z,\tau_nX)\simeq\ext^1_{\aa}(\Omega^{i-1}Z,\tau_nX)\stackrel{\XBBC(3)}{\simeq}D\underline{\bb}(\Omega^{n-1}X,\Omega^{i-1}Z)\stackrel{\XBA(1)}{\simeq}D\underline{\bb}(\Omega^{n-i}X,Z)$, which is $\stackrel{\XBA(2)}{\simeq}D\ext^{n-i}_{\aa}(X,Z)$ if $n>i$. The isomorphisms for $\tau^-_n$ are given dually.\rule{5pt}{10pt}

\vskip.5em{\bf\XBCB\ Corollary }{\it
$\tau_n$ and $\tau_n^-$ give mutually quasi-inverse equivalences $\tau_{n}:\underline{\xx}_{n-1}\to\overline{\yy}_{n-1}$ and $\tau_{n}^-:\overline{\yy}_{n-1}\to\underline{\xx}_{n-1}$.}

\vskip.5em{\sc Proof }
For any $X\in\underline{\xx}_{n-1}$, $\ext^{i}_{\aa}(\ii,\tau_nX)\stackrel{\XBCA}{\simeq}D\ext^{n-i}_{\aa}(X,\ii)=0$ holds for any $i$ ($0<i<n$). Thus $\tau_n$ gives a functor $\underline{\xx}_{n-1}\to\overline{\yy}_{n-1}$, which is full and faithful by \XBA(1) and \XBBC(3). Dually, $\tau^-_n$ gives a full and faithful functor $\underline{\yy}_{n-1}\to\overline{\xx}_{n-1}$. Since $\underline{\bb}(X,\ )\stackrel{\XBCA}{\simeq}D\ext^n_{\aa}(\ ,\tau_nX)\stackrel{\XBCA}{\simeq}\underline{\bb}(\tau_n^-\circ\tau_nX,\ )$ holds, $\tau_n^-\circ\tau_n$ is isomorphic to the identity functor. Dually, $\tau_n\circ\tau_n^-$ is isomorphic to the identity functor.\rule{5pt}{10pt}

\vskip.5em{\bf\XBD\ }
Let $\cc$ be a functorially finite subcategory of $\bb$, and $l\ge0$. We call $\cc$ an {\it $l$-orthogonal} subcategory of $\bb$ if $\cc\perp_l\cc$ holds, and a {\it maximal $l$-orthogonal} subcategory of $\bb$ if $\cc=\cc^{\perp_l}\cap\bb={}^{\perp_l}\cc\cap\bb$ holds. We call $M\in\bb$ {\it maximal $l$-orthogonal} (resp. {\it $l$-orthogonal}) if so is $\add M$. Any maximal $l$-orthogonal subcategory $\cc$ of $\bb$ satisfies $\pp\cup\ii\subseteq\cc\subseteq\xx_l\cap\yy_l$. Since $\underline{\cc}$ and $\overline{\cc}$ are functorially finite subcategories of $\underline{\bb}$ and $\overline{\bb}$ respectively, $\underline{\cc}$ and $\overline{\cc}$ are dualizing $R$-varieties by \XAB. Thus $\mod\underline{\cc}$, $\mod\underline{\cc}^{op}$, $\mod\overline{\cc}$ and $\mod\overline{\cc}^{op}$ are closed under kernels, cokernels and extensions by \XAA. We have the following characterizations of maximal $l$-orthogonal subcategories [I7;2.2.2].

\vskip.5em{\bf\XBDA\ Proposition }{\it
Let $\cc$ be a functorially finite subcategory of $\bb$. Then the conditions (1), (2-$i$) and (3-$i$) are equivalent for any $i$ ($0\le i\le l$).

(1) $\cc$ is a maximal $l$-orthogonal subcategory of $\bb$.

(2-$0$) $\resdim{\cc}{\bb}\le l$, $\cc\perp_l\cc$ and $\pp\cup\ii\subseteq\cc$.

(2-$i$) $\resdim{\cc}{(\cc^{\perp_i}\cap\bb)}\le l-i$, $\cc\perp_l\cc$ and $\pp\cup\ii\subseteq\cc$.

(2-$l$) $\cc=\cc^{\perp_l}\cap\bb$ and $\pp\subseteq\cc$.

(3-$0$) $\resdim{\cc^{op}}{\bb^{op}}\le l$, $\cc\perp_l\cc$ and $\pp\cup\ii\subseteq\cc$.

(3-$i$) $\resdim{\cc^{op}}{({}^{\perp_i}\cc\cap\bb)^{op}}\le l-i$, $\cc\perp_l\cc$ and $\pp\cup\ii\subseteq\cc$.

(3-$l$) $\cc={}^{\perp_l}\cc\cap\bb$ and $\ii\subseteq\cc$.}


\vskip.5em{\bf\XBE\ }In the rest of this section, let $\cc$ be a maximal $(n-1)$-orthogonal subcategory of $\bb$ ($n\ge1$). Assume that $\cc$ is {\it Krull-Schmidt}, i.e. any object of $\cc$ is isomorphic to a finite direct sum of objects whose endomorphism rings are local.

\vskip.5em{\bf\XBEA\ }The following fundamental theorem follows from previous results in \XBC\ and \XBD\ (cf. [I7;2.3,2.3.1,2.2.3;3.5.2]).

\vskip.5em{\bf Theorem }{\it
(1) {\rm ($n$-Auslander-Reiten translation)} For any $X\in\cc$, $\tau_nX\in\cc$ and $\tau_n^-X\in\cc$ hold. Thus $\tau_n:\underline{\cc}\to\overline{\cc}$ and $\tau_n^-:\overline{\cc}\to\underline{\cc}$ are mutually quasi-inverse equivalences.

(2) {\rm ($n$-Auslander-Reiten duality)} There exist functorial isomorphisms $\overline{\cc}(Y,\tau_nX)\simeq D\ext^n_{\aa}(X,Y)\simeq\underline{\cc}(\tau^-_nY,X)$ for any $X,Y\in\cc$.

(3) $\resdim{\cc}{\bb}\le n-1$ and $\resdim{\cc^{op}}{\bb^{op}}\le n-1$ hold.

(4) $X\mapsto\ext^n_{\aa}(\ ,X)$ gives an equivalence $\overline{\cc}\to\ii(\mod\underline{\cc})$, and $X\mapsto\ext^n_{\aa}(X,\ )$ gives an equivalence $\underline{\cc}^{op}\to\ii(\mod\overline{\cc}^{op})$.}


\vskip.5em{\bf\XBEB\ }For any $F\in\mod\cc$, take a projective resolution $\cc(\ ,Y)\stackrel{\cdot f}{\to}\cc(\ ,X)\to F\to0$. Define $\alpha F\in\mod\cc^{op}$ by the exact sequence $0\to\alpha F\to\cc(X,\ )\stackrel{f\cdot}{\to}\cc(Y,\ )$. Then $\alpha$ gives a left exact functor $\alpha:\mod\cc\rightarrow\mod\cc^{op}$. Define $\alpha:\mod\cc^{op}\rightarrow\mod\cc$ dually. We denote by $\der{n}:\mod\cc\leftrightarrow\mod\cc^{op}$ the $n$-th derived functor of $\alpha$ [FGR]. Then we have the following theorem (see [I7;3.6.1]). 

\vskip.5em{\bf Theorem }{\it
Let $\cc$ be a maximal $(n-1)$-orthogonal subcategory of $\bb$ ($n\ge1$).

(1) Any $0\neq F\in\mod\underline{\cc}$ satisfies $\pd{}_{\cc}F=n+1$ and $\der{i}F=0$ ($i\neq n+1$). Any $0\neq G\in\mod\overline{\cc}^{op}$ satisfies $\pd{}_{\cc^{op}}G=n+1$ and $\der{i}G=0$ ($i\neq n+1$).

(2) $\der{n+1}$ gives a duality $\mod\underline{\cc}\leftrightarrow\mod\overline{\cc}^{op}$, and the equivalence $D\der{n+1}:\mod\underline{\cc}\leftrightarrow\mod\overline{\cc}$ coincides with the equivalence induced by $\tau_n:\underline{\cc}\to\overline{\cc}$.}

\vskip.5em{\bf\XBEC\ }
We can show the following theorem (see [I7;3.3,3.3.1]).

\vskip.5em{\bf Theorem }($n$-almost split sequence) {\it
Let $\cc$ be a maximal $(n-1)$-orthogonal subcategory of $\bb$ ($n\ge1$). Fix any non-projective $X\in\ind\cc$ (resp. non-injective $Y\in\ind\cc$).

(1) There exists an exact sequence ${\bf A}:0\to Y\stackrel{f_n}{\to}C_{n-1}\stackrel{f_{n-1}}{\to}\cdots\stackrel{f_1}{\to}C_0\stackrel{f_0}{\to}X\to0$ with terms in $\cc$ such that $f_i\in J_{\cc}$ and the following sequences are exact.
{\small\begin{eqnarray*}
0\to\cc(\ ,Y)\stackrel{\cdot f_n}{\to}\cc(\ ,C_{n-1})\stackrel{\cdot f_{n-1}}{\to}\cdots\stackrel{\cdot f_1}{\to}\cc(\ ,C_0)\stackrel{\cdot f_0}{\to}J_{\cc}(\ ,X)\to0\\
0\to\cc(X,\ )\stackrel{f_0\cdot}{\to}\cc(C_0,\ )\stackrel{f_1\cdot}{\to}\cdots\stackrel{f_{n-1}\cdot}{\to}\cc(C_{n-1},\ )\stackrel{f_n\cdot}{\to}J_{\cc}(Y,\ )\to0
\end{eqnarray*}}
Such ${\bf A}$ is unique up to isomorphisms of complexes, and satisfies $Y\simeq\tau_nX$ and $X\simeq\tau_n^-Y$. 

(2) The simple modules $F:=\cc/J_{\cc}(\ ,X)$ and $G:=\cc/J_{\cc}(Y,\ )$ satisfy $\pd{}_{\cc}F=n+1=\pd_{\cc^{op}}G$, $\der{i}F=0=\der{i}G$ ($i\neq n+1$),  $F=\der{n+1}G$ and $G=\der{n+1}F$.}

\vskip.5em{\bf\XBF\ }In the rest of this section, we fix $m\ge0$ and impose the conditions below on $\bb$.

\vskip.5em{\bf Proposition }{\it
For $m\ge0$, the following conditions for $\bb$ are equivalent.

(1) $\Omega^mX\in\bb$ holds for any $X\in\aa$.

(2) $\bb$ is a contravariantly finite subcategory of $\aa$ with $\resdim{\bb}{\aa}\le m$.

(3) If $0\to Y\to B_{m-1}\to\cdots\to B_0\to X\to0$ is an exact sequence in $\aa$ with $B_i\in\bb$, then $Y\in\bb$ holds.}

\vskip.5em{\sc Proof }
(3)$\Rightarrow$(1) is obvious. We will show (2)$\Rightarrow$(3). Take the following commutative diagram of exact sequences, where the upper sequence is a right $\bb$-resolution of $X$.
\[\begin{diag}
0&\RA{}&C_m&\RA{}&C_{m-1}&\RA{}&C_{m-2}&\RA{}&&\RA{}&C_0&\RA{}&X&\RA{}&0\\
&&\uparrow&&\uparrow^{}&&\uparrow^{}&&\cdots&&\uparrow^{}&&\parallel\\
0&\RA{}&Y&\RA{}&B_{m-1}&\RA{}&B_{m-2}&\RA{}&&\RA{}&B_0&\RA{}&X&\RA{}&0
\end{diag}\]

Taking mapping cone, we have an exact sequence $0\to Y\to C_m\oplus B_{m-1}\to\cdots\to C_1\oplus B_0\to C_0\to0$. Since $\bb$ is a resolving subcategory of $\aa$, we obtain $Y\in\bb$.

(1)$\Rightarrow$(2) Fix any $X\in\aa$. By Auslander-Buchweitz approximation theory [AB;1.1], there exists an exact sequence $0\to I_m\to\cdots\to I_1\to B_0\to X\to0$ with $I_i\in\ii$ and $B_0\in\bb$. It is easily checked that this is a right $\bb$-resolution of $X$.\rule{5pt}{10pt}

\vskip.5em{\bf\XBFA\ }We can show the following theorem by a similar argument to [I7;3.6.2].

\vskip.5em{\bf Theorem }{\it
Any maximal $(n-1)$-orthogonal subcategory $\cc$ ($n\ge1$) of $\bb$ satisfies $\gl(\mod\cc)\le\max\{n+1,m\}$.}

\vskip.5em{\bf\XBFB\ }We end this section by pointing out the interesting result below, which realizes the category $\cc^{\perp_{i}}\cap\bb$ as the category of syzygies.

\vskip.5em{\bf Theorem }{\it
Let $\cc=\add M$ be an $(n-1)$-orthogonal subcategory of $\bb$. Assume $\pp\cup\ii\subseteq\cc$ and that $\Gamma:=\endm_{\cc}(M)$ is a noetherian ring. For any $i$ ($0\le i\le n-1$), we have full and faithful functors $\fff:=\bb(M,\ ):\cc^{\perp_{i}}\cap\bb\to\Omega^{i+2}(\mod\Gamma)$ and $\ggg:=\bb(\ ,M):{}^{\perp_{i}}\cc\cap\bb\to\Omega^{i+2}(\mod\Gamma^{op})$ such that $\fff=(\ )^*\circ\ggg$ and $\ggg=(\ )^*\circ\fff$ for $(\ )^*=\hom_\Gamma(\ ,\Gamma)$. If $m-2\le i$, then $\fff$ and $\ggg$ are equivalences.}

\vskip.5em{\sc Proof }
We only show the assertion for $\fff$. For any $X\in\bb$, take a right $\cc$-resolution $C_1\to C_0\to X\to0$, which is exact by $\pp\subseteq\cc$. We have exact sequences $0\to\bb(X,\ )\to\bb(C_0,\ )\to\bb(C_1,\ )$ and $0\to{}_\Gamma(\fff X,\fff(\ ))\to{}_\Gamma(\fff C_0,\fff(\ ))\to{}_\Gamma(\fff C_1,\fff(\ ))$ on $\bb$. Since $\bb(C_i,\ )={}_\Gamma(\fff C_i,\fff(\ ))$ holds on $\bb$, $\fff$ is full and faithful and $\ggg=(\ )^*\circ\fff$ holds. For any $X\in\cc^{\perp_{i}}\cap\bb$, take an injective resolution ${\bf I}:0\to X\to I_0\to\cdots\to I_{i+1}$ in $\bb$. Since $\cc\perp_{i}X$ holds, $\fff{\bf I}:0\to\fff X\to\fff I_0\to\cdots\to\fff I_{i+1}$ is an exact sequence with $\fff I_j\in\add{}_\Gamma\Gamma$. Thus $\fff X\in\Omega^{i+2}(\mod\Gamma)$ holds. 

Assume $m-2\le i$. For any $Y\in\Omega^{i+2}(\mod\Gamma)$, take an exact sequence ${\bf P}:0\to Y\to P_{i+1}\to\cdots\to P_0$ with $P_j\in\add{}_\Gamma\Gamma$. Since $\fff$ gives an equivalence $\cc\to\add{}_\Gamma\Gamma$, we can take a complex ${\bf C}:C_{i+1}\stackrel{f_{i+1}}{\to}\cdots\stackrel{f_1}{\to}C_0$ with $C_j\in\cc$ such that $\fff{\bf C}$ is isomorphic to ${\bf P}$. Since $\pp\subseteq\cc$ and ${\bf P}$ is exact, ${\bf C}$ is also exact. Put $X_j:=\Ker f_{j-1}$. Inductively, we can easily show that $\cc\perp_{j}X_{j+2}$ holds for any $j$ by using exactness of $\fff{\bf C}$. In particular, $X_{i+2}\in\cc^{\perp_{i}}\cap\bb$ holds by $m-2\le i$, and $\fff X_{i+2}=Y$ holds.\rule{5pt}{10pt}

\vskip.5em{\bf\XC\ Orders, cotilting modules and Auslander-type conditions }

The aim of this section is to define a pair $(\aa,\bb)$ to which we will apply our theory in \S2, and to give preliminary facts which we will use in preceding sections.

\vskip.5em{\bf\XCA\ }
Throughout this section, let $R$ be a complete regular local ring of dimension $d$ and $\Lambda$ a module-finite $R$-algebra. We call $\Lambda$ an {\it isolated singularity} [A3] if $\gl\Lambda\otimes_RR_{\dn{p}}=\height\dn{p}$ holds for any non-maximal prime ideal $\dn{p}$ of $R$. We call a left $\Lambda$-module $M$ {\it Cohen-Macaulay} if it is a projective $R$-module. We denote by $\cm\Lambda$ the category of Cohen-Macaulay $\Lambda$-modules. Then $D_d:=\hom_R(\ ,R)$ gives a duality $\cm\Lambda\leftrightarrow\cm\Lambda^{op}$. We call $\Lambda$ an {\it $R$-order} (or {\it Cohen-Macaulay $R$-algebra}) if $\Lambda\in\cm\Lambda$ [A2,3]. A typical example of an order is a commutative complete local Cohen-Macaulay ring $\Lambda$ containing a field since such $\Lambda$ contains a complete regular local subring $R$ [Ma;29.4]. Let ${\bf E}:0\to R\to E_0\to\cdots\to E_d\to0$ be a minimal injective resolution of the $R$-module $R$. We denote by $D:=\hom_R(\ ,E_d)$ the Matlis dual. Put $(\ )^*:=\hom_\Lambda(\ ,\Lambda)$ and denote by $\nu_\Lambda:=D_d\circ(\ )^*$ the Nakayama functor and by $\nu^-_\Lambda:=(\ )^*\circ D_d$ the inverse Nakayama functor. If $\Lambda$ is an $R$-order, then $\nu_\Lambda$ and $\nu^-_\Lambda$ give mutually inverse equivalences $\add{}_\Lambda\Lambda\leftrightarrow\add{}_\Lambda(D_d\Lambda)$. The following observation in [I7;2.5.1] is useful.

\vskip.5em{\bf\XCAA\ Proposition }{\it
Let $\Lambda$ be an $R$-order which is an isolated singularity, $X,Y\in\cm\Lambda$ and $2\le n\le d$. Then $\depth_R\hom_\Lambda(X,Y)\ge n$ if and only if $X\perp_{n-2}Y$.}

\vskip.5em{\bf\XCAB\ }Let $\Lambda$ be a module-finite $R$-algebra which is an isolated singularity and $M\in\cm\Lambda$. Let us recall the method of Goto and Nishida [GN1] to construct a minimal injective resolution of $M$ in $\mod\Lambda$ from a minimal projective resolution ${\bf P}:\cdots\to P_1\to P_0\to D_dM\to0$ of $D_dM$. We have exact sequences $M\otimes_R{\bf E}:0\to M\to M\otimes_RE_0\to\cdots\to M\otimes_RE_{d-1}\to M\otimes_RE_{d}\to0$ and $D{\bf P}:0\to M\otimes_RE_{d}\to DP_0\to DP_1\to\cdots$. Connecting them, we obtain a minimal injective resolution
\[0\to M\to M\otimes_RE_0\to\cdots\to M\otimes_RE_{d-1}\to DP_0\to DP_1\to\cdots\]
of $M$ in $\mod\Lambda$. Thus $\id{}_\Lambda M=\pd(D_dM)_\Lambda+d$ holds. In particular, if $\gl\Lambda=d$, then $\cm\Lambda\subseteq\add{}_\Lambda\Lambda$ holds.

\vskip.5em{\bf\XCB\ }Let $\Lambda$ be an $R$-order. For $m\ge d$, we call $T\in\cm\Lambda$ an {\it $m$-cotilting module} [AR3][N] if $T\perp T$ (\XBA) and there exist exact sequences $0\to T\to I_0\to\cdots\to I_{m-d}\to0$ and $0\to T_{m-d}\to\cdots\to T_0\to D_d\Lambda\to0$ with $I_i\in\add{}_\Lambda(D_d\Lambda)$ and $T_i\in\add{}_\Lambda T$. It is easy to check the facts below by \XCAB\ and \XCAA.

(1) $\id{}_\Lambda T\le m$ and ${}^\perp T\subseteq\cm\Lambda$ hold, and $\endm_\Lambda(T)$ is an $R$-order.

(2) $T\in\cm\Lambda$ is a $d$-cotilting module if and only if $\add{}_\Lambda T=\add{}_\Lambda (D_d\Lambda)$.

\vskip.5em{\bf\XCBA\ }Let us recall the following classical {\it cotilting theorem} [M][Ha].

\vskip.5em{\bf Proposition }{\it
Let $T$ be an $m$-cotilting $\Lambda$-module and ${\Lambda^\prime}:=\endm_\Lambda(T)^{op}$. Then $T$ is an $m$-cotilting ${\Lambda^\prime}$-module. We have mutually quasi-inverse equivalences $\hom_\Lambda(\ ,T):{}^\perp({}_\Lambda T)\to{}^\perp({}_{\Lambda^\prime} T)$ and $\hom_{{\Lambda^\prime}}(\ ,T):{}^\perp({}_{\Lambda^\prime} T)\to{}^\perp({}_\Lambda T)$ which preserve $\ext^i$ for any $i\ge0$.}

\vskip.5em{\bf\XCBB\ }We can apply our theory in \S\XB\ to $(\aa,\bb):=(\mod\Lambda,{}^\perp T)$ by the following version of a theorem of Auslander-Reiten [AR3].

\vskip.5em{\bf Proposition }{\it Let $\Lambda$ be an $R$-order which is an isolated singularity, $T$ an $m$-cotilting $\Lambda$-module, $\aa:=\mod\Lambda$ and $\bb:={}^\perp T$. Then the following assertions hold.

(1) $\bb$ is a enough injective resolving subcategory of $\aa$ with $\ii(\bb)=\add T$.

(2) $\bb$ is a functorially finite subcategory of $\aa$ with $\resdim{\bb}{\aa}\le m$.

(3) $\underline{\bb}$ and $\overline{\bb}$ are dualizing $R$-varieties.}

\vskip.5em{\sc Proof }
(1) It is easily checked that $\bb$ is resolving with $T\in\ii$. For any $X\in\bb$, take an injection $X\stackrel{a}{\to}(D_d\Lambda)^l$ in $\cm\Lambda$ by \XCB(1). Take an exact sequence $0\to T_{m-d}^l\to\cdots\to T_0^l\stackrel{b}{\to}(D_d\Lambda)^l\to0$ in \XCB. Then $a$ factors through $b$ by $X\perp T$. Thus $X$ is a submodule of $T_0^l$. We can take an exact sequence $0\to X\stackrel{c}{\to}T^\prime\to Y\to0$ such that $c$ is a left $(\add T)$-resolution of $X$. Applying ${}_\Lambda(\ ,T)$, we obtain $Y\in{}^\perp T=\bb$. Thus $\bb$ is enough injectives with $\ii=\add T$.

(2) Since $\Omega^m\aa\subseteq\bb$ holds by $\id{}_\Lambda T\le m$, $\bb$ is a contravariantly finite subcategory of $\aa$ with $\resdim{\bb}{\aa}\le m$ by \XBF. We will show that $\bb$ is a covariantly finite subcategory of $\aa$. Put ${\Lambda^\prime}:=\endm_\Lambda(T)^{op}$ and $\bb^\prime:={}^\perp({}_{\Lambda^\prime} T)\subseteq\aa^\prime:=\mod{\Lambda^\prime}$. Since $T$ is an $m$-cotilting ${\Lambda^\prime}$-module by \XCBA, $\bb^\prime$ is a contravariantly finite subcategory of $\aa^\prime$. Fix $X\in\aa$. Let $B^\prime\stackrel{a}{\to}{}_\Lambda(X,T)$ be a right $\bb^\prime$-resolution. It is easily checked that the composition $X\to{}_{{\Lambda^\prime}}({}_\Lambda(X,T),T)\stackrel{a\cdot}{\to}{}_{{\Lambda^\prime}}(B^\prime,T)$ is a left $\bb$-resolution of $X$.

(3) Put $\bb_0:=\cm\Lambda$. Since $\Lambda$ is an isolated singularity, it is well-known that $\bb_0$ satisfies the conditions in \XBBC\ (e.g. [A2;8.7][AR5;2.4]). Since $\bb$ is a functorially finite subcategory of $\bb_0$ by (2), it is easily checked that $\underline{\bb}$ is that of $\underline{\bb}_0$. Thus $\underline{\bb}$ is a dualizing $R$-variety by \XAB, and so is $\overline{\bb}$ by \XBBC.\rule{5pt}{10pt}

\vskip.5em{\bf\XCBC\ }Let us recall the theorem [I7;3.4.4] below, which tells us that higher dimensional Auslander-Reiten theory for the case $d=m=n+1$ is quite peculiar. It means that $\cc$ has sequences which have properties like $n$-almost split sequences and connect projective modules and injective modules. We notice that $f_0$ below is not surjective in general.

\vskip.5em{\bf Theorem }($n$-fundamental sequence) {\it
Let $\bb$ be in \XCBB\ and $\cc$ a maximal $(n-1)$-orthogonal subcategory of $\bb$. Assume $d=m=n+1$. Fix any $X\in\cc$ (resp. $Y\in\cc$).

(1) There exists an exact sequence ${\bf A}:0\to Y\stackrel{f_n}{\to}C_{n-1}\stackrel{f_{n-1}}{\to}\cdots\stackrel{f_1}{\to}C_0\stackrel{f_0}{\to}X$ with terms in $\cc$ such that $f_i\in J_{\cc}$ and the following sequences are exact.
{\small\begin{eqnarray*}
0\to\cc(\ ,Y)\stackrel{\cdot f_n}{\to}\cc(\ ,C_{n-1})\stackrel{\cdot f_{n-1}}{\to}\cdots\stackrel{\cdot f_1}{\to}\cc(\ ,C_0)\stackrel{\cdot f_0}{\to}J_{\cc}(\ ,X)\to0\\
0\to\cc(X,\ )\stackrel{f_0\cdot}{\to}\cc(C_0,\ )\stackrel{f_1\cdot}{\to}\cdots\stackrel{f_{n-1}\cdot}{\to}\cc(C_{n-1},\ )\stackrel{f_n\cdot}{\to}J_{\cc}(Y,\ )\to0
\end{eqnarray*}}
Such ${\bf A}$ is unique up to isomorphisms of complexes, and satisfies $Y\simeq\nu_\Lambda X$ and $X\simeq\nu_\Lambda ^-Y$. 

(2) The simple modules $F:=\cc/J_{\cc}(\ ,X)$ and $G:=\cc/J_{\cc}(Y,\ )$ satisfy $\pd{}_{\cc}F=n+1=\pd_{\cc^{op}}G$, $\der{i}F=0=\der{i}G$ ($i\neq n+1$), $F=\der{n+1}G$ and $G=\der{n+1}F$.}

\vskip.5em{\bf\XCC\ Definition }Let us introduce certain Auslander-type conditions on selfinjective resolutions, which will play a crucial role in this paper (see \XDB). Let $\Gamma$ be a noetherian ring and $0\to\Gamma\to I_0\to I_1\to\cdots$ a minimal injective resolution of the $\Gamma$-module $\Gamma$. We say that $\Gamma$ satisfies the {\it $(m,n)$-condition} if $\fd{}_\Gamma I_i<m$ holds for any $i$ ($i<n$) [I2,4]. We can state many well-known homological conditions in terms of our $(m,n)$-conditions. For example, the {\it dominant dimension} $\domdim\Gamma:=\inf\{ i\ge0\ |\ \fd{}_\Gamma I_i\neq0\}$ [T][H] of $\Gamma$ is the maximal number $n$ such that $\Gamma$ satisfies the $(1,n)$-condition. Moreover, recall that $\Gamma$ is called {\it $n$-Gorenstein} if $\fd{}_\Gamma I_i\le i$ holds for any $i$ ($0\le i<n$) [FGR][B][AR4][C]. This is equivalent to that $\Gamma$ satisfies the $(i,i)$-condition for any $i$ ($0<i\le n$). We notice that our $(m,n)$-condition itself is not left-right symmetric. We say that $\Gamma$ satisfies the {\it two-sided $(m,n)$-condition} if $\Gamma$ and $\Gamma^{op}$ satisfies the $(m,n)$-condition.

\vskip.5em{\bf\XCCA\ Proposition }{\it
Let $\Gamma$ be an $R$-order which is an isolated singularity, and $0\to\Gamma\to I_0\to I_1\to\cdots$ a minimal injective resolution in $\cm\Gamma$.

(1) $\Gamma$ is $d$-Gorenstein.

(2) $\Gamma$ satisfies the $(m,n)$-condition if and only if $\pd{}_\Gamma I_i<m-d$ for any $i$ ($i<n-d$).

(3) If $I\in\add(D_d\Gamma)$ satisfies $\pd{}_\Gamma I\le n$, then $I\in\add(\bigoplus_{i=0}^{n}I_i)$ holds.}

\vskip.5em{\sc Proof }
(1) and (2) follow by \XCAB\ since $\fd{}_\Gamma(M\otimes_RE_i)=i$ ($i<d$) and $\fd{}_\Gamma(I_i\otimes_RE_d)=\pd{}_\Gamma I_i+d$ ($i\ge0$) hold by [GN1]. Miyachi's theorem [Mi] implies (3).\rule{5pt}{10pt}

\vskip.5em{\bf\XCD\ }Let us introduce {\it $m$-extension pairs}, which will be used in \XDD. As we will see in \XCDC, they are closely related to $m$-cotilting modules.

\vskip.5em{\bf\XCDA\ Proposition }{\it
Let $\Gamma$ be a module-finite $R$-algebra, and $e$ and $f$ idempotents of $\Gamma$ such that $\Gamma f\in\cm\Gamma$ and $e\Gamma\in\cm\Gamma^{op}$. Put $P:=\Gamma f$ and $I:=D_d(e\Gamma)$. The conditions (1) and (2) below are equivalent.

(1) Put $\underline{\Gamma}:=\Gamma/\Gamma e\Gamma$ and $\overline{\Gamma}:=\Gamma/\Gamma f\Gamma$. For any $i\ge0$, $\ext^i_\Gamma(\ ,\Gamma)$ gives functors $\mod\underline{\Gamma}\to\mod\overline{\Gamma}^{op}$ and $\mod\overline{\Gamma}^{op}\to\mod\underline{\Gamma}$.

(2) There exist exact sequences $0\to P\to I_0\to I_1\to\cdots$ and $\cdots\to P_1\to P_{0}\to I\to 0$ with $I_i\in\add I$ and $P_i\in\add P$.

If the conditions above and $\id{}_\Gamma P\le m$ and $\id{}_\Gamma(D_dI)\le m$ are satisfied, we call $(P,I)$ an {\rm $m$-extension pair}. Then we can assume $I_{m-d+1}=0=P_{m-d+1}$ in (2).}

\vskip.5em{\sc Proof }
Notice that $Y\in\mod\Gamma^{op}$ is contained in $\mod\overline{\Gamma}^{op}$ if and only if $Y\otimes_\Gamma P=0$.

(1)$\Rightarrow$(2) Let $\cdots\to P_1^\prime\to P_0^\prime\to D_dP\to0$ be a minimal projective resolution. By \XCAB, we have a minimal injective resolution $0\to P\to P\otimes_RE_0\to\cdots\to P\otimes_RE_{d-1}\to DP_0^\prime\to DP_1^\prime\to\cdots$. Fix any simple $S\in\mod\underline{\Gamma}$. Since $\ext^i_\Gamma(S,\Gamma)\in\mod\overline{\Gamma}^{op}$ holds for any $i$, we obtain $D(P_i^\prime\otimes_\Gamma S)={}_\Gamma(S,DP_i^\prime)=\ext^{i+d}_\Gamma(S,P)=\ext^{i+d}_\Gamma(S,\Gamma)\otimes_\Gamma P=0$ by $\ext^{i+d}_\Gamma(S,\Gamma)\in\mod\overline{\Gamma}^{op}$. Thus $P_i^\prime\otimes_\Gamma S=0$ holds, and $P_i^\prime\in\add(e\Gamma)_\Gamma$ for any $i$.

(2)$\Rightarrow$(1) By \XCDB\ below, $\ext^j_\Gamma(X,I_i)=0$ holds for any $X\in\mod\underline{\Gamma}$ and $i,j\ge0$. Since we have an exact sequence $0\to P\to I_0\to I_1\to\cdots$, we have $\ext^j_\Gamma(X,\Gamma)\otimes_\Gamma P=\ext^j_\Gamma(X,P)=0$ for any $j$. Thus $\ext^j_\Gamma(X,\Gamma)\in\mod\overline{\Gamma}^{op}$.\rule{5pt}{10pt}

\vskip.5em{\bf\XCDB\ Lemma }{\it
Let $\Gamma$ be a module-finite $R$-algebra, and $e$ an idempotent of $\Gamma$ such that $e\Gamma\in\cm\Gamma^{op}$. Put $\underline{\Gamma}:=\Gamma/\Gamma e\Gamma$ and $I:=D_d(e\Gamma)$. Then $\ext^i_\Gamma(X,I\otimes_RY)=0$ holds for any $i\ge0$, $X\in\mod\underline{\Gamma}$ and $Y\in\Mod R$.}

\vskip.5em{\sc Proof }
Put $Q:=\Gamma e$ and $\qqq:=\hom_\Gamma(Q,\ )$. We have a functorial isomorphism ${}_\Gamma(\ ,I\otimes_RY)={}_R(Q^*\otimes_\Gamma\ ,Y)={}_R(\qqq(\ ),Y)$. Let ${\bf A}:\cdots\to P_1\to P_0\to X\to0$ be a projective resolution of $X\in\mod\underline{\Gamma}$. We have an exact sequence $\qqq{\bf A}:\cdots\to\qqq P_{1}\to\qqq P_0\to0$. Since $\qqq P_i\in\cm R$ holds for any $i$, the complex $\qqq{\bf A}$ splits as a complex of $R$-modules. Thus we obtain an exact sequence ${}_R(\qqq{\bf A},Y):0\to{}_R(\qqq P_{0},Y)\to{}_R(\qqq P_1,Y)\to\cdots$. By the remark above, ${}_\Gamma({\bf A},I\otimes_RY):0\to{}_\Gamma(P_0,I\otimes_RY)\to{}_\Gamma(P_1,I\otimes_RY)\to\cdots$ is exact. Thus $\ext^i_\Gamma(X,I\otimes_RY)=0$.\rule{5pt}{10pt}

\vskip.5em{\bf\XCDC\ Proposition }{\it
(1) Let $\Lambda$ be an $R$-order, $T$ an $m$-cotilting $\Lambda$-module and $\Lambda\oplus T\in\add{}_\Lambda M\subseteq{}^\perp T$. Put $\Gamma:=\endm_\Lambda(M)$, $P:=\hom_\Lambda(M,T)$ and $I:=D_dM$. Then $(P,I)$ is an $m$-extension pair of $\Gamma$-modules.

(2) Let $\Gamma$ be an module-finite $R$-algebra and $(P,I)$ an $m$-extension pair. Put $Q:=\nu^-_\Gamma I$, $\Lambda:=\endm_\Gamma(Q)$, $M:=\hom_\Gamma(Q,\Gamma)=D_dI$ and $T:=\hom_\Gamma(Q,P)$. Then $\Lambda$ is an $R$-order, $T$ is an $m$-cotilting $\Lambda$-module and $\Lambda\oplus T\in\add{}_\Lambda M\subseteq{}^\perp T$.}

\vskip.5em{\sc Proof }
(1) Take exact sequences ${\bf I}:0\to T\to I_0\to\cdots\to I_{m-d}\to0$ and ${\bf T}:0\to T_{m-d}\to\cdots\to T_0\to D_d\Lambda\to0$ in \XCB. Since $M\in{}^\perp T$, we obtain exact sequences ${}_\Lambda(M,{\bf I}):0\to P\to {}_\Lambda(M,I_0)\to\cdots\to {}_\Lambda(M,I_{m-d})\to0$ with ${}_\Lambda(M,I_i)\in\add{}_\Gamma I$ and ${}_\Lambda(M,{\bf T}):0\to {}_\Lambda(M,T_{m-d})\to\cdots\to {}_\Lambda(M,T_0)\to I\to0$ with ${}_\Lambda(M,T_i)\in\add{}_\Gamma P$. Thus $(P,I)$ is an $m$-extension pair.

(2) By our assumption, $M=D_dI\in\cm R$ holds and $\Lambda$ is an $R$-order. Put $\qqq:=\hom_\Gamma(Q,\ )$. Then $\qqq I=D_d\Lambda$ and ${}_\Lambda(M,\qqq I)=D_dM=I$ hold. Take exact sequences ${\bf I}:0\to P\to I_0\to\cdots\to I_{m-d}\to0$ and ${\bf P}:0\to P_{m-d}\to\cdots\to P_{0}\to I\to 0$ with $I_i\in\add{}_\Gamma I$ and $P_i\in\add{}_\Gamma P$. We have an exact sequence $\qqq{\bf I}:0\to T\to\qqq I_0\to\cdots\to\qqq I_{m-d}\to0$ with $\qqq I_i\in\add{}_\Lambda(D_d\Lambda)$, which gives an injective resolution of $T$. Thus $\id{}_\Lambda T\le m$ holds. Since ${}_\Lambda(M,\qqq{\bf I})$ is isomorphic to the exact sequence ${\bf I}$ by the remark above, $M\in{}^\perp T$ holds. On the other hand, we have an exact sequence $\qqq{\bf P}:0\to\qqq P_{m-d}\to\cdots\to\qqq P_0\to D_d\Lambda\to 0$ with $\qqq P_i\in\add{}_\Lambda T$. Thus $T$ is an $m$-cotilting $\Lambda$-module. Since $Q\oplus P\in\add{}_\Gamma\Gamma$, we obtain $\Lambda\oplus T\in\add{}_\Lambda M$.\rule{5pt}{10pt}

\vskip.5em{\bf\XCE\ }Let us introduce {\it $n$-superprojective} modules, which will be used in \XDD.

\vskip.5em{\bf\XCEA\ Proposition }{\it
Let $\Gamma$ be a module-finite $R$-algebra, and $e$ an idempotent of $\Gamma$ such that $e\Gamma\in\cm\Gamma^{op}$. Put $Q:=\Gamma e$, $I:=D_d(e\Gamma)$ and $\underline{\Gamma}:=\Gamma/\Gamma e\Gamma$. For $n\ge 1$, the conditions (1)--(3) below are equivalent.

(1) $\grade{{}_\Gamma X}\ge n+1$ holds for any $X\in\mod\underline{\Gamma}$.

(2) There exists an exact sequence $0\to\Gamma\to I_0\to\cdots\to I_n$ with $I_i\in\add{}_\Gamma I$.

(3) Put $\Lambda:=\endm_\Gamma(Q)$. Then the functor $\qqq:=\hom_\Gamma(Q,\ ):\add{}_\Gamma\Gamma\to\mod\Lambda$ is full and faithful and $\qqq\Gamma\in\mod\Lambda$ is $(n-1)$-orthogonal.

If the conditions above are satisfied, we call $Q$ {\rm $n$-superprojective}. Moreover, if $n\ge d$ and $\Gamma$ is an isolated singularity, then the condition (4) below is also equivalent.

(4) $\Gamma$ is an $R$-order with an injective resolution $0\to\Gamma\to I_0\to\cdots\to I_{n-d}$ in $\cm\Gamma$ with $I_i\in\add{}_\Gamma I$.}

\vskip.5em{\bf\XCEB\ }For the proof, we need the following easy lemma.

\vskip.5em{\bf Lemma }{\it
Let ${\bf P}:P_{n+1}\to P_n\to\cdots\to P_0\to0$ be a complex with $P_i\in\add{}_\Gamma\Gamma$, and $H_i$ the homology of ${\bf P}$ at $P_i$. If $\grade{{}_\Gamma H_i}>n-i$ holds for any $i$ ($0\le i\le n$), then ${\bf P}^*:0\to P_{0}^*\to\cdots\to P_{n}^*\to P_{n+1}^*$ is exact for $(\ )^*=\hom_\Gamma(\ ,\Gamma)$.}

\vskip.5em{\bf\XCEC\ Proof of \XCEA\ }
By our assumption, $\Lambda$ is an $R$-order and $\qqq$ gives a functor $\qqq=\hom_\Gamma(Q,\ ):\add{}_\Gamma\Gamma\to\cm\Lambda$.

(2)$\Rightarrow$(1) By \XCDB, $\ext^j_\Gamma(X,I_i)=0$ holds for any $i,j\ge0$. Since we have an exact sequence $0\to\Gamma\to I_0\to\cdots\to I_n$, $\ext^j_\Gamma(X,\Gamma)=0$ holds for any $j$ ($j\le n$).

(1)$\Rightarrow$(3) Let ${\bf Q}:Q_n\to\cdots\to Q_0\to \Gamma$ be a right $(\add{}_\Gamma Q)$-resolution of $\Gamma$ and $H_i$ the homology of ${\bf Q}$ at $Q_i$. Since $\qqq H_i=0$ holds, we obtain $\grade{{}_\Gamma H_i}>n$ for any $i$. By \XCEB, ${\bf Q}^*:0\to\Gamma\to Q_0^*\to\cdots\to Q_n^*$ is exact. On the other hand, we have a projective resolution $\qqq{\bf Q}:\qqq Q_n\to\cdots\to\qqq Q_0\to\qqq\Gamma\to0$ of $\qqq\Gamma\in\mod\Lambda$. Thus we have an exact sequence ${}_\Lambda(\qqq{\bf Q},\qqq\Gamma):0\to{}_\Lambda(\qqq\Gamma,\qqq\Gamma)\to{}_\Lambda(\qqq Q_0,\qqq\Gamma)\to\cdots\to{}_\Lambda(\qqq Q_n,\qqq\Gamma)$ with a homology $\ext^i_\Lambda(\qqq\Gamma,\qqq\Gamma)$ at ${}_\Lambda(\qqq Q_i,\qqq\Gamma)$ for any $i>0$. Since we have the following commutative diagram of complexes, $\qqq$ is full and faithful and $\qqq\Gamma$ is $(n-1)$-orthogonal.
\[\begin{diag}
{\bf Q}^*:&0&\longrightarrow&\Gamma&\longrightarrow&Q_0^*&\longrightarrow&Q_1^*&\longrightarrow&\cdots&\longrightarrow&Q_n^*\\
&&&\downarrow&&\parallel&&\parallel&&&&\parallel\\
{}_\Lambda(\qqq{\bf Q},\qqq\Gamma):&0&\longrightarrow&{}_\Lambda(\qqq\Gamma,\qqq\Gamma)&\longrightarrow&{}_\Lambda(\qqq Q_0,\qqq\Gamma)&\longrightarrow&{}_\Lambda(\qqq Q_1,\qqq\Gamma)&\longrightarrow&\cdots&\longrightarrow&{}_\Lambda(\qqq Q_n,\qqq\Gamma)
\end{diag}\]

(3)$\Rightarrow$(2) Since $\Lambda$ is an $R$-order, we can take an injective resolution ${\bf A}:0\to\qqq\Gamma\to I_0^\prime\to\cdots\to I_n^\prime$ in $\cm\Lambda$. Since $\qqq\Gamma$ is $(n-1)$-orthogonal, ${}_\Lambda(\qqq\Gamma,{\bf A})$ is exact with ${}_\Lambda(\qqq\Gamma,I_i^\prime)\in\add{}_\Gamma(D_d\qqq\Gamma)=\add{}_\Gamma I$. Thus ${}_\Lambda(\qqq\Gamma,{\bf A})$ gives the desired sequence.

We will show the assertion for (4). Obviously (2)$\Rightarrow$(4) holds. We will show (4)$\Rightarrow$(1). Let ${\bf E}:0\to R\stackrel{f_0}{\to}E_0\stackrel{f_1}{\to}\cdots\stackrel{f_d}{\to}E_d\to0$ be a minimal injective resolution of the $R$-module $R$ and $F_i:=\Cok f_{i-1}$ for $i\le d$. By \XCDB, $\ext^j_\Gamma(X,I\otimes_RF_i)=0$ holds for any $j\ge0$. Since we have an exact sequence $0\to\Gamma\otimes_RF_i\to I_0\otimes_RF_i\to\cdots\to I_{n-d}\otimes_RF_i$ with $I_i\in\add{}_\Gamma I$, $\ext^j_\Gamma(X,\Gamma\otimes_RF_i)=0$ holds for any $j$ ($j\le n-d$). Since the $i$-th cosyzygy of $\Gamma$ in $\mod\Gamma$ is $\Gamma\otimes_RF_i$ by \XCAB, we obtain $\ext^{i+j}_\Gamma(X,\Gamma)=\ext^j_\Gamma(X,\Gamma\otimes_RF_i)=0$ for any $i$ ($i\le d$) and $j$ ($j\le n-d$). Thus $\grade{{}_\Gamma X}>n$ holds.\rule{5pt}{10pt}

\vskip.5em{\bf\XD\ Higher dimensional Auslander algebras }

Throughout this section, fix a complete regular local ring $R$ of dimension $d\ge0$. 

\vskip.5em{\bf\XDA\ Definition }
Let $m\ge d$ and $n\ge 1$. An {\it Auslander} (resp. {\it quasi-Auslander}) {\it triple of type $(d,m,n)$} is a triple $(\Lambda,M,T)$ which satisfies (1)--(3) (resp. (1)(2) and (3)$^\prime$) below.

(1) $\Lambda$ is an $R$-order which is an isolated singularity, and $T,M\in\mod\Lambda$.

(2) $T$ is an $m$-cotilting $\Lambda$-module.

(3) $\add{}_\Lambda M$ is a maximal $(n-1)$-orthogonal subcategory of ${}^\perp T$.

(3)$^\prime$ $\add{}_\Lambda M$ is an $(n-1)$-orthogonal subcategory of ${}^\perp T$ and contains $\Lambda$ and $T$.

We call an $R$-algebra $\Gamma$ an {\it Auslander} (resp. {\it quasi-Auslander}) {\it algebra of type $(d,m,n)$} if there exists an Auslander (resp. quasi-Auslander) triple $(\Lambda,M,T)$ of type $(d,m,n)$ such that $\Gamma=\endm_\Lambda(M)$.

\vskip.5em{\bf\XDAA\ }
(1) We will consider triples $(\Lambda,M_1,M_2)$ of a noetherian ring $\Lambda$ and $M_i\in\mod\Lambda$. We say that two triples $(\Lambda^i,M_1^i,M_2^i)$ ($i=1,2$) are {\it equivalent} if there exists an equivalence $\mod\Lambda^1\to\mod\Lambda^2$ which induces equivalences $\add{}_{\Lambda^1}M_j^1\to\add{}_{\Lambda^2}M_j^2$ for $j=1,2$.

(2) For $m\ge d$ and $n\ge 1$, we denote by $\dn{A}_{m,n}$ (resp. $\dn{A}_{m,n}^q$) the set of equivalence classes of Auslander (resp. quasi-Auslander) triples of type $(d,m,n)$. Then $\dn{A}_{m,n}^q\supseteq\dn{A}_{m,n}\supseteq\dn{A}_{m^\prime,n}$ and $\dn{A}_{m,n}^q\supseteq\dn{A}_{m^\prime,n^\prime}^q$ hold for any $m\ge m^\prime$ and $n\le n^\prime$. For any element of $\dn{A}_{m,n}$ (resp. $\dn{A}_{m,n}^q$), its associated Auslander (resp. quasi-Auslander) algebra is uniquely determined up to Morita-equivalence.

\vskip.5em{\bf\XDB\ Main Results }
We collect our main results which will be proved in \XDF. 

\vskip.5em{\bf\XDBA\ }For the case $m\le n$, we can give the homological characterization of (quasi-) Auslander algebras of type $(d,m,n)$ below by using Auslander-type condition in \XCC. The case $(d,m,n)=(0,0,1)$ is given by Auslander [A1] and Auslander-Solberg [ASo], the case $(d,m,n)=(1,1,1)$ is given by Auslander-Roggenkamp [ARo], and the case $(d,m,n)=(0,1,1)$ is given by the author [I5].

\vskip.5em{\bf Theorem }{\it
Let $\Gamma$ be an $R$-algebra. If $m\le n$, then $\Gamma$ is an Auslander (resp. quasi-Auslander) algebra of type $(d,m,n)$ if and only if $\Gamma$ is an $R$-order which is an isolated singularity and satisfies the two-sided $(m+1,n+1)$-condition and $\gl\Gamma\le n+1$ (resp. the two-sided $(m+1,n+1)$-condition).}

\vskip.5em{\bf\XDBB\ }A more explicit result is given by Auslander correspondence below for the case $m\le n$. A more general result for arbitrary case will be given in \XDDA.

\vskip.5em{\bf Theorem }(Auslander correspondence of type $(d,m,n)$) {\it
Assume $m\le n$. Then the map $(\Lambda,M,T)\mapsto\endm_\Lambda(M)$ gives a bijection from $\dn{A}_{m,n}$ (resp. $\dn{A}_{m,n}^q$) to the set of Morita-equivalence classes of $R$-orders $\Gamma$ which are isolated singularities and satisfy the two-sided $(m+1,n+1)$-condition and $\gl\Gamma\le n+1$ (resp. the two-sided $(m+1,n+1)$-condition). In particular, two triples $(\Lambda_i,M_i,T_i)\in\dn{A}_{m,n}^q$ ($i=1,2$) are equivalent if and only if two categories $\add{}_{\Lambda_i}M_i$ ($i=1,2$) are equivalent.}

\vskip.5em{\bf\XDBC\ }
Let us study the case $d=m>n$. Then $\add{}_\Lambda T=\add{}_\Lambda(D_d\Lambda)$ and ${}^\perp T=\cm\Lambda$ hold by \XCB(2). In this case, we can give the homological characterization of Auslander algebras below. In particular, putting $n:=1$, we obtain an answer to M. Artin's question [Ar] to give a homological characterization of endomorphism rings $\endm_\Lambda(M)$ of additive generators $M$ of $\cm\Lambda$ for representation-finite orders $\Lambda$.

\vskip.5em{\bf Theorem }{\it
Let $\Gamma$ be an $R$-algebra. If $d>n$, then $\Gamma$ is an Auslander algebra of type $(d,d,n)$ if and only if the conditions (1) and (2) below hold.

(1) $\Gamma$ is a module-finite $R$-algebra which is an isolated singularity with $\gl\Gamma=d$ and $\depth{}_R\Gamma\ge n+1$.

(2) There exists an idempotent $e$ of $\Gamma$ such that $e\Gamma\in\cm\Gamma^{op}$, $\underline{\Gamma}:=\Gamma/\Gamma e\Gamma$ is artinian, and $\pd{}_\Gamma X\le n+1$ holds for any $X\in\mod\underline{\Gamma}$.}

\vskip.5em{\bf\XDBD\ }For the arbitrary case, we can give a homological characterization of Auslander algebras in the theorem below. These conditions strongly reflect properties of maximal $(n-1)$-orthogonal subcategories studied in \S\XB.

\vskip.5em{\bf Theorem }{\it
Let $\Gamma$ be an $R$-algebra, $m\ge d$ and $n\ge1$. Then $\Gamma$ is an Auslander algebra of type $(d,m,n)$ if and only if the conditions (1)--(4) below hold.

(1) $\Gamma$ is a module-finite $R$-algebra which is an isolated singularity with $\gl\Gamma\le\max\{n+1,m\}$.

(2) There exist idempotents $e$ and $f$ of $\Gamma$ such that $e\Gamma\in\cm\Gamma^{op}$, $\id(e\Gamma)_\Gamma\le m$, $\Gamma f\in\cm\Gamma$ and $\id{}_\Gamma(\Gamma f)\le m$.

(3) Put $\underline{\Gamma}:=\Gamma/\Gamma e\Gamma$ and $\overline{\Gamma}:=\Gamma/\Gamma f\Gamma$. Any $0\neq X\in\mod\underline{\Gamma}$ satisfies $\pd{}_\Gamma X=\grade{{}_\Gamma X}=n+1$, and any $0\neq Y\in\mod\overline{\Gamma}^{op}$ satisfies $\pd Y_\Gamma=\grade{Y_\Gamma}=n+1$.

(4) $\ext^{n+1}_\Gamma(\ ,\Gamma)$ gives a duality $\mod\underline{\Gamma}\leftrightarrow\mod\overline{\Gamma}^{op}$.}

\vskip.5em{\bf\XDC\ }Let us start with collecting properties of (quasi-)Auslander algebras.

(1) If $(\Lambda,M,T)\in\dn{A}_{m,n}$ (resp. $\dn{A}_{m,n}^q$), then $(\endm_\Lambda(T)^{op},\hom_\Lambda(M,T),T)\in\dn{A}_{m,n}$ (resp. $\dn{A}_{m,n}^q$) and $\endm_{\endm_\Lambda(T)^{op}}(\hom_\Lambda(M,T))=\endm_\Lambda(M)^{op}$ hold by \XCBA. Consequently, $\Gamma$ is an Auslander (resp. quasi-Auslander) algebra of type $(d,m,n)$ if and only if so is $\Gamma^{op}$.

(2) Assume $m:=\gl\Lambda<\infty$. Then $\Lambda$ is an $m$-cotilting $\Lambda$-module with ${}^\perp\Lambda=\add\Lambda$. Since $(\Lambda,\Lambda,\Lambda)\in\dn{A}_{m,n}$ holds for any $n\ge1$, $\Lambda$ is an Auslander algebra of type $(d,m,n)$. We call the equivalence class of such a triple {\it trivial}.

\vskip.5em{\bf\XDCA\ Proposition }{\it
Let $(\Lambda,M,T)\in\dn{A}_{m,n}^q$ and $\Gamma:=\endm_\Lambda(M)$. Then (1)--(6) below hold. If $(\Lambda,M,T)\in\dn{A}_{m,n}$, then (7)--(9) below hold.

(1) $\Gamma$ is a module-finite $R$-algebra which is an isolated singularity.

(2) We have mutually inverse equivalences $\mmm:=\hom_\Lambda(M,\ ):\add{}_\Lambda M\to\add{}_\Gamma\Gamma$ and $\qqq:=\hom_\Gamma(Q,\ ):\add{}_\Gamma\Gamma\to\add{}_\Lambda M$ for $Q:=\mmm\Lambda$.

(3) $M\in\cm\Gamma^{op}\cap\add\Gamma_\Gamma$ and $P:=\hom_\Lambda(M,T)\in\cm\Gamma\cap\add{}_\Gamma\Gamma$. For $I:=D_dM$, $(P,I)$ is an $m$-extension pair (\XCDA) and $Q:=\nu^-_\Gamma I(=\mmm\Lambda)$ is $n$-superprojective (\XCEA). 

(4) A left $(\add{}_\Gamma P)$-resolution $0\to\Gamma\to P_0\to\cdots\to P_n$ of ${}_\Gamma\Gamma$ and a left $(\add(D_dI)_\Gamma)$-resolution $0\to\Gamma\to P_0^\prime\to\cdots\to P_n^\prime$ of $\Gamma_\Gamma$ are exact.

(5) $\Gamma$ satisfies $\depth_R\Gamma\ge\min\{n+1,d\}$ and the two-sided $(m+1,n+1)$-condition.

(6) Take idempotents $e$ and $f$ of $\Gamma$ such that $\add{}_\Gamma I=\add{}_\Gamma D_d(e\Gamma)$ and $\add{}_\Gamma P=\add{}_\Gamma(\Gamma f)$. Then $\underline{\Gamma}:=\Gamma/\Gamma e\Gamma$ and $\overline{\Gamma}:=\Gamma/\Gamma f\Gamma$ are artin algebras with the following commutative diagram for $\cc:=\add M$, $\underline{\cc}:=\cc/[\add{}_\Lambda\Lambda]$ and $\overline{\cc}:=\cc/[\add{}_\Lambda T]$.
\[\begin{diag}
\mod\underline{\Gamma}&\ \ \ \subset\ \ \ &\mod\Gamma&\ \ \ \stackrel{\ext^i_\Gamma(\ ,\Gamma)}{\longleftrightarrow}\ \ \ &\mod\Gamma^{op}&\ \ \ \supset\ \ \ &\mod\overline{\Gamma}^{op}\\
\uparrow{\scriptstyle\wr}&&\uparrow{\scriptstyle\wr}&&\uparrow{\scriptstyle\wr}&&\uparrow{\scriptstyle\wr}\\
\mod\underline{\cc}&\ \ \ \subset\ \ \ &\mod\cc&\ \ \ \stackrel{\der{i}}{\longleftrightarrow}\ \ \ &\mod\cc^{op}&\ \ \ \supset\ \ \ &\mod\overline{\cc}^{op}\end{diag}\]

\vskip-.5em
(7) $d\le\gl\Gamma\le\max\{n+1,\id{}_\Lambda T\}$ holds, and the right equality holds if $(\Lambda,M,T)$ is non-trivial (\XDC(2)).

(8) Any $0\neq X\in\mod\underline{\Gamma}$ satisfies $\pd{}_\Gamma X=\grade{{}_\Gamma X}=n+1$, and any $0\neq Y\in\mod\overline{\Gamma}^{op}$ satisfies $\pd Y_\Gamma=\grade{Y_\Gamma}=n+1$.

(9) $\ext^{n+1}_\Gamma(\ ,\Gamma)$ gives a duality $\mod\underline{\Gamma}\leftrightarrow\mod\overline{\Gamma}^{op}$.}

\vskip.5em{\sc Proof }
(1) Obviously $\Gamma$ is a finitely generated $R$-module. For any non-maximal prime ideal $\dn{p}$ of $R$, $M_{\dn{p}}$ is a progenerator of $\Lambda_{\dn{p}}$ (e.g. [R;3.5]). Thus $\Gamma_{\dn{p}}=\endm_{\Lambda_{\dn{p}}}(M_{\dn{p}})$ is Morita-equivalent to $\Lambda_{\dn{p}}$. This implies that $\Gamma$ is an isolated singularity.

(2) Obviously $\mmm$ is an equivalence. Moreover, $\qqq\circ\mmm={}_\Gamma(Q,{}_\Lambda(M,\ ))={}_\Lambda(M\otimes_\Gamma Q,\ )={}_\Lambda(M\otimes_\Gamma\hom_\Lambda(M,\Lambda),\ )={}_\Lambda(\Lambda,\ )=1$ holds.

(3) $(P,I)$ is an $m$-extension pair by \XCDC(1), and $Q$ is $n$-superprojective by \XCEA(3).

(4) Take exact sequences ${\bf T}:0\to M\to T_0\to\cdots\to T_n$ and ${\bf P}:P_n\to\cdots\to P_0\to M\to0$ with $T_i\in\add T$ and $P_i\in\add\Lambda$. Then ${}_\Lambda(M,{\bf P})$ and ${}_\Lambda({\bf P},M)$ are exact by $M\perp_{n-1}M$. They are left $(\add{}_\Gamma P)$ and $(\add(D_dI)_\Gamma)$-resolutions by (2).

(5) The former assertion follows by (4). Take an exact sequence $0\to\Gamma\to I_0\to\cdots\to I_n$ with $I_i\in\add I$ in \XCEA(2). By \XCAB, $0\to I\to I\otimes_RE_0\to\cdots\to I\otimes_RE_d\to0$ gives a minimal injective resolution of $I$. Since $\pd{}_\Gamma I\le m-d$ holds, $\fd{}_\Gamma I\otimes_RE_i\le m$ holds. The mapping cone gives an injective resolution of $\Gamma$, which shows that $\Gamma$ satisfies the $(m+1,n+1)$-condition. By \XDC(1), $\Gamma^{op}$ satisfies the $(m+1,n+1)$-condition.

(6) We have an equivalence $\mod\cc\to\mod\Gamma$ given by $F\mapsto F(M)$, which makes our diagram commutative. Since $\Lambda$ is an isolated singularity, $\underline{\Gamma}$ and $\overline{\Gamma}$ are artin algebras.

(7) Since $\gl\Lambda\ge d$ holds by [R;3.2], the former assertion follows by \XBFA. We will show the latter one. Since $(\Lambda,M,T)$ is non-trivial, there exists non-projective $X\in\ind\cc$. Then the $\Gamma$-module $F:=\cc/J_{\cc}(M,X)$ satisfies $\pd{}_\Gamma F=n+1$ by \XBEB. Let $0\to C_l\to\cdots\to C_1\to J_\Lambda\to0$ be a minimal right $\cc$-resolution of $J_\Lambda$ with $C_l\neq0$. Then the $\Gamma$-module $G:=\cc/J_{\cc}(M,\Lambda)$ satisfies $\pd{}_\Gamma G=l$. Since we have an exact sequence $0\to C_l\to\cdots\to C_1\to\Lambda\to\Lambda/J_\Lambda\to0$ with $C_i\oplus\Lambda\perp T$, we have $l\ge d$ and $\ext_\Lambda^{l+1}(\Lambda/J_\Lambda,T)=0$. By \XCAB, $\id{}_\Lambda T\le l$ holds. Thus $\gl\Gamma\ge\max\{n+1,\id{}_\Lambda T\}$.

(8) follows by \XBEB(1), and (9) follows by \XBEB(2).\rule{5pt}{10pt}

\vskip.5em{\bf\XDD\ Definition }
To prove the theorems in \XDB, it is convenient to introduce $\dn{B}_{m,n}$ as follows. We denote by $\dn{B}_{m,n}$ (resp. $\dn{B}_{m,n}^q$) the set of equivalence classes (\XDAA(1)) of triples $(\Gamma,P,I)$ which satisfies the conditions (1)--(3) (resp. (1) and (2)) below.

(1) $\Gamma$ is a module-finite $R$-algebra which is an isolated singularity.

(2) $(P,I)$ is an $m$-extension pair (\XCDA) and $Q:=\nu^-_\Gamma I$ is $n$-superprojective (\XCEA).

(3) $\gl\Gamma\le\max\{n+1,m\}$ and any $X\in\mod\underline{\Gamma}$ (\XCDA) satisfies $\pd{}_\Gamma X\le n+1$.

We will show in \XDDD\ that the sets $\dn{B}_{m,n}$ and $\dn{B}_{m,n}^q$ are `left-right symmetric'.

\vskip.5em{\bf\XDDA\ Theorem }(Auslander correspondence of type $(d,m,n)$){\it

(1) There exists a bijection $\alpha:\dn{A}_{m,n}^q\to\dn{B}_{m,n}^q$ for any $m\ge d$ and $n\ge1$. It is given by $\alpha(\Lambda,M,T):=(\endm_\Lambda(M),\hom_\Lambda(M,T),D_dM)$, and the converse is given by $\alpha^{-1}(\Gamma,P,I):=(\endm_\Gamma(Q),\hom_\Gamma(Q,\Gamma),\hom_\Gamma(Q,P))$ for $Q:=\nu^-_\Gamma I$.

(2) $\alpha$ gives a bijection $\alpha:\dn{A}_{m,n}\to\dn{B}_{m,n}$ for any $m\ge d$ and $n\ge1$.

(3) $\Gamma$ is an Auslander (resp. quasi-Auslander) algebra of type $(d,m,n)$ if and only if $(\Gamma,P,I)\in\dn{B}_{m,n}$ (resp. $\dn{B}_{m,n}^q$) holds for some $(P,I)$.
\[\begin{array}{ccccccccc}
&\Lambda,&&M,&&T,&&D_d\Lambda&\\
{}^{\hom_\Lambda(M,\ )}&\downarrow\uparrow&&\downarrow\uparrow&&\downarrow\uparrow&&\downarrow\uparrow&{}^{\hom_\Gamma(Q,\ )}\\
&Q,&&\Gamma,&&P,&&I&
\end{array}\]}

\vskip-.5em{\bf\XDDB\ Lemma }{\it
(1) For any $(\Lambda,M,T)\in\dn{A}_{m,n}^q$, put $\Gamma:=\endm_\Lambda(M)$, $P:=\hom_\Lambda(M,T)$, $I:=D_dM$ and $Q:=\hom_\Lambda(M,\Lambda)$. Then $(\Gamma,P,I)\in\dn{B}_{m,n}^q$. Moreover, $\Lambda=\endm_\Gamma(Q)$, $M=\hom_\Gamma(Q,\Gamma)$ and $T=\hom_\Gamma(Q,P)$ hold.

(2) For any $(\Gamma,P,I)\in\dn{B}_{m,n}^q$, put $Q:=\nu^-_\Gamma I\in\add{}_\Gamma\Gamma$, $\Lambda:=\endm_\Gamma(Q)$, $M:=\hom_\Gamma(Q,\Gamma)=D_dI$ and $T:=\hom_\Gamma(Q,P)$. Then $(\Lambda,M,T)\in\dn{A}_{m,n}^q$. Moreover, $\Gamma=\endm_\Lambda(M)$, $P=\hom_\Lambda(M,T)$ and $I=\hom_\Lambda(M,D_d\Lambda)$ hold.}

\vskip.5em{\sc Proof }
(1) $(P,I)$ is an $m$-extension pair by \XCDC(1). The latter assertion follows by \XDCA(2). Thus $\qqq:\add{}_\Gamma\Gamma\to\mod\Lambda$ is full and faithful and $M=\qqq\Gamma$ is $(n-1)$-orthogonal. Hence $Q$ is $n$-superprojective, and $(\Gamma,P,I)\in\dn{B}_{m,n}^q$ holds.

(2) By \XCEA(3), $\qqq:\add{}_\Gamma\Gamma\to\mod\Lambda$ is full and faithful and $M=\qqq\Gamma$ is $(n-1)$-orthogonal. Thus the former assertion holds by \XCDC(2). Since ${}_\Lambda(M,\qqq(\ ))={}_\Lambda(\qqq\Gamma,\qqq(\ ))={}_\Gamma(\Gamma,\ )=1$ holds on $\add{}_\Gamma\Gamma$, the latter assertion follows.\rule{5pt}{10pt}

\vskip.5em{\bf\XDDC\ Proof of \XDDA\ }
(1) follows by \XDDB. We will show (2). Fix $(\Lambda,M,T)\in\dn{A}_{m,n}^q$ and the corresponding $(\Gamma,P,I)\in\dn{B}_{m,n}^q$. If $(\Lambda,M,T)\in\dn{A}_{m,n}$, then \XDD(3) holds by \XDCA(7)(8), so $(\Gamma,P,I)\in\dn{B}_{m,n}$ holds. We will show that $(\Gamma,P,I)\in\dn{B}_{m,n}$ implies $(\Lambda,M,T)\in\dn{A}_{m,n}$, i.e. $\cc:=\add{}_\Lambda M$ is a maximal $(n-1)$-orthogonal subcategory of ${}^\perp T$. By \XBDA, we only have to show that any $X\in\cc^{\perp_{n-1}}\cap{}^\perp T$ satisfies $X\in\cc$. Put $g:=\max\{n+1,m\}$ and $\mmm:={}_\Lambda(M,\ )$. Take an exact sequence ${\bf T}:0\to X\stackrel{f_0}{\to}T_0\stackrel{f_1}{\to}\cdots\stackrel{f_{g-1}}{\to}T_{g-1}$ with $T_i\in\add{}_\Lambda T$. Applying $\mmm$, we obtain a complex $\mmm{\bf T}:0\to\mmm X\stackrel{}{\to}\mmm T_0\stackrel{}{\to}\cdots\stackrel{}{\to}\mmm T_{g-1}$ with $\mmm T_i\in\add{}_\Gamma P$. Let $H_i$ be the homology of $\mmm{\bf T}$ at $\mmm T_i$. Since ${}_\Gamma(Q,\mmm(\ ))$ is the identity functor, ${}_\Gamma(Q,H_i)=0$ holds for any $i$. By \XDD(3), $\pd{}_\Gamma H_i\le n+1$ holds for any $i$. Put $X_0:=X$ and $X_i:=\Cok f_{i-1}$ for $i>0$. Inductively, we will show that $\pd{}_\Gamma\mmm X_i\le i$ holds for any $i$ ($n-1\le i\le g$). This is true for $i=g$ by $\gl\Gamma\le g$. Assume that $\pd{}_\Gamma\mmm X_i\le i$ holds for some $i$ ($n\le i\le g$). We have an exact sequence $0\to\mmm X_{i-1}\to\mmm T_{i-1}\to\mmm X_i\stackrel{g_i}{\to}H_i\to 0$. Since $\pd{}_\Gamma\mmm X_{i}\le i$ and $\pd{}_\Gamma H_i\le n+1$ hold, we obtain $\pd{}_\Gamma\Ker g_i\le i$. Thus $\mmm T_{i-1}\in\add{}_\Gamma\Gamma$ implies $\pd{}_\Gamma\mmm X_{i-1}\le i-1$. In particular, $\pd{}_\Gamma\mmm X_{n-1}\le n-1$ holds. Since $M$ is $(n-1)$-orthogonal, $H_i=0$ for any $i$ ($0\le i<n$). Thus $\mmm X$ is an $(n-1)$-st syzygy of $\mmm X_{n-1}$, and $\mmm X\in\add{}_\Gamma\Gamma$ holds. Thus we obtain $X={}_\Gamma(Q,\mmm X)\in\cc$.\rule{5pt}{10pt}

\vskip.5em{\bf\XDDD\ Corollary }{\it
Let $m\ge d$ and $n\ge1$. If $(\Gamma,P,I)\in\dn{B}_{m,n}$ (resp. $\dn{B}_{m,n}^q$), then $(\Gamma^{op},D_dI,D_dP)\in\dn{B}_{m,n}$ (resp. $\dn{B}_{m,n}^q$).}

\vskip.5em{\sc Proof }
Put $(\Lambda,M,T):=\alpha^{-1}(\Gamma,P,I)\in\dn{A}_{m,n}$ (resp. $\dn{A}_{m,n}^q$) and $\Lambda^\prime:=\endm_\Lambda(T)^{op}$. Then $(\Lambda^\prime,P,T)=(\endm_\Lambda(T)^{op},{}_\Lambda(M,T),T)\in\dn{A}_{m,n}$ (resp. $\dn{A}_{m,n}^q$) and $\endm_{\Lambda^\prime}(P)=\Gamma^{op}$ hold by \XDC(1). Since ${}_{\Lambda^\prime}(P,T)=M=D_dI$ holds by \XCBA, we obtain $(\Gamma,D_dI,D_dP)=\alpha(\Lambda^\prime,P,T)\in\dn{B}_{m,n}$ (resp. $\dn{B}_{m,n}^q$).\rule{5pt}{10pt}

\vskip.5em{\bf\XDE\ }The following proposition connects \XDBB\ and \XDDA.

\vskip.5em{\bf Proposition }{\it
If $m\le n$, then the map $(\Gamma,P,I)\mapsto\Gamma$ gives a bijection from $\dn{B}_{m,n}$ (resp. $\dn{B}_{m,n}^q$) to the set of Morita-equivalence classes of $R$-orders $\Gamma$ which are isolated singularities and satisfy the two-sided $(m+1,n+1)$-condition and $\gl\Gamma\le n+1$ (resp. the two-sided $(m+1,n+1)$-condition).}

\vskip.5em{\sc Proof }
We only have to show the assertion for $\dn{B}_{m,n}^q$.

(i) Fix $(\Gamma,P,I)\in\dn{B}_{m,n}^q$. Applying \XDCA(5) to $(\Lambda,M,T):=\alpha^{-1}(\Gamma,P,I)\in\dn{A}_{m,n}^q$, $\Gamma$ is an $R$-order and satisfies the two-sided $(m+1,n+1)$-condition.

(ii) Fix an $R$-order $\Gamma$ satisfying the two-sided $(m+1,n+1)$-condition. Take minimal injective resolutions $0\to\Gamma\to I_0\to\cdots\to I_{n-d}$ in $\cm\Gamma$ and $0\to\Gamma\to J_0\to\cdots\to J_{n-d}$ in $\cm\Gamma^{op}$. Put $I:=\bigoplus_{i=0}^{n-d}I_i$, $Q:=\nu_\Gamma^-I$, $J:=\bigoplus_{i=0}^{n-d}J_i$ and $P:=D_dJ$. By \XCEA(4), $Q$ is $n$-superprojective. Since $\Gamma$ satisfies the two-sided $(m+1,n+1)$-condition, \XCCA(2) implies $\pd{}_\Gamma I\le m-d$ and $\pd J_\Gamma\le m-d$. Thus $(P,I)$ is an $m$-extension pair by $m\le n$. Thus $(\Gamma,P,I)\in\dn{B}_{m,n}^q$.

We will show that any $(\Gamma,P^\prime,I^\prime)\in\dn{B}_{m,n}^q$ is equivalent to $(\Gamma,P,I)$. Since there exists an injective resolution $0\to\Gamma\to I_0^\prime\to\cdots\to I_{n-d}^\prime$ in $\cm\Gamma$ with $I_i^\prime\in\add{}_\Gamma I^\prime$ by \XCEA(4), $I\in\add{}_\Gamma I^\prime$ holds. Since $\pd{}_\Gamma I^\prime\le m-d\le n-d$ holds by \XCDA, $I^\prime\in\add{}_\Gamma I$ holds by \XCCA(3). Thus $\add{}_\Gamma I^\prime=\add{}_\Gamma I$. A dual argument shows $\add{}_\Gamma P^\prime=\add{}_\Gamma P$.\rule{5pt}{10pt}

\vskip.5em{\bf\XDF\ Proof of Main results }
\XDBB\ follows immediately by \XDDA\ and \XDE, and \XDBA\ follows by \XDBB. We will show \XDBD. The `only if' part follows by \XDCA, and `if' part follows by \XDDA(3) since $(\Gamma,\Gamma f,D_d(e\Gamma))\in\dn{B}_{m,n}$ holds by \XCDA(1) and \XCEA(1).

We will show \XDBC. The `only if' part follows by \XDCA. We will show the `if' part. Put $I:=D_d(e\Gamma)$. Since $\depth{}_RI=d$ and $\gl\Gamma=d$ hold by (2), $I\in\add{}_\Gamma\Gamma$ holds by \XCAB. Take an idempotent $f$ of $\Gamma$ such that $\add{}_\Gamma I=\add{}_\Gamma P$ holds for $P:=\Gamma f$. Then $(P,I)$ is a $d$-extension pair by $\id{}_\Gamma P=d$ and $\id(D_dI)_\Gamma=d$. Since any $X\in\mod\underline{\Gamma}$ has finite length by (2), we obtain $\grade{{}_\Gamma X}\ge n+1$ by $\depth{}_R\Gamma\ge n+1$ in (1). Thus $Q$ is $n$-superprojective, and $(\Gamma,P,I)\in\dn{B}_{m,n}$ holds. The assertion follows by \XDDA(3).\rule{5pt}{10pt}

\vskip.5em{\bf\XDG\ }Recall that higher dimensional Auslander-Reiten theory for the case $d=m=n+1$ is quite peculiar (\XCBC). Correspondingly, Auslander algebras of type $(d,d,d-1)$ have a very nice homological characterizations below. In particular, the condition (3) below means that $\Gamma$ is a (non-local and non-graded version of) {\it Artin-Schelter regular ring of dimension $d$} [ArS]. The symmetry of projective resolutions of simple modules over Artin-Schelter regular rings corresponds to the selfduality of $(d-1)$-almost split sequences and $(d-1)$-fundamental sequences. See \XEBB\ and \XFA\ for examples. The equivalence of (1) and (2) for $d=2$ is a theorem of Auslander [Ar][RV].

\vskip.5em{\bf Theorem }{\it
For a module-finite $R$-algebra $\Gamma$, the conditions below are equivalent.

(1) $\Gamma$ is an Auslander algebra of type $(d,d,d-1)$.

(2) $\Gamma$ is an $R$-order with $\gl\Gamma=d$.

(3) $\gl\Gamma=d$ holds, and any simple $\Gamma$-module $S$ satisfies $\ext^{i}_\Gamma(S,\Gamma)=0$ ($i\neq d$) and $\ext^d_\Gamma(S,\Gamma)$ is a simple $\Gamma^{op}$-module.

(4) Opposite side version of (3).}

\vskip.5em{\sc Proof }
(1)$\Rightarrow$(3) Immediate from \XDCA(7) and \XCBC.

(3)$\Rightarrow$(2) Immediate from $\depth_R\Gamma=\inf\{i\ge0\ |\ \ext^i_\Gamma(\Gamma/J_\Gamma,\Gamma)\neq0\}$ [GN2;3.2].

(2)$\Rightarrow$(1) Since $\cm\Gamma=\add{}_\Gamma\Gamma$ holds by \XCAB, $\Gamma\in\cm\Gamma$ is maximal $(d-1)$-orthogonal.\rule{5pt}{10pt}

\vskip.5em{\bf\XDGA\ }Let us recall the proposition below [I4;6.3]. An important example of such $\Gamma$ is an Auslander algebra of type $(d,n,n)$, which satisfies (1) and (2) below by \XDBA\ and \XBEC\ respectively. In this sense, the two-sided $(n+1,n+1)$-condition means the existence of $n$-almost split sequences homologically.

\vskip.5em{\bf Proposition }{\it 
For a noetherian ring $\Gamma$ with $\gl\Gamma=n+1$, the conditions below are equivalent.

(1) $\Gamma$ satisfies the two-sided $(n+1,n+1)$-condition.

(2) Any simple $\Gamma$-module (resp. $\Gamma^{op}$-module) $S$ with $\pd S=n+1$ satisfies $\ext^{i}_\Gamma(S,\Gamma)=0$ ($i\neq n+1$) and $\ext^{n+1}_\Gamma(S,\Gamma)$ is a simple $\Gamma^{op}$-module (resp. $\Gamma$-module).}

\vskip.5em{\bf\XDH\ }Let us study how to get all quasi-Auslander triples with a fixed quasi-Auslander algebra. For any automorphism $\phi\in\aut(\Lambda)$, we have the induced auto-equivalence $\phi:\mod\Lambda\to\mod\Lambda$. Then any quasi-Auslander triple $(\Lambda,M,T)$ gives another quasi-Auslander triple $(\Lambda,\phi(M),\phi(T))$ with the same quasi-Auslander algebra. Since $\phi(X)$ is isomorphic to $X$ for any $\phi\in\inn(\Lambda)$ and $X\in\mod\Lambda$, this action of $\aut(\Lambda)$ factors through $\out(\Lambda)$. By the theorem below, $\out(\Lambda)$ is sufficient for our purpose. We call a triple $(\Lambda,M,T)$ {\it basic} if all algebras $\Lambda$, $\endm_\Lambda(M)$ and $\endm_\Lambda(T)$ are basic.

\vskip.5em{\bf Theorem }{\it
Let $(\Lambda,M_i,T_i)$ be a basic quasi-Auslander triple of type $(d,m,n)$ ($i=1,2$). Assume $m\le n$. Then $\endm_\Lambda(M_1)$ is isomorphic to $\endm_\Lambda(M_2)$ if and only if there exists $\phi\in\out(\Lambda)$ such that $\phi(M_1)$ and $\phi(T_1)$ are isomorphic to $M_2$ and $T_2$ respectively.}

\vskip.5em{\sc Proof }By \XDBB\ and our definition of $\dn{A}_{m,n}^q$ in \XDAA, there exists an auto-equivalence $\fff:\mod\Lambda\to\mod\Lambda$ which induces equivalences $\add M_1\to\add M_2$ and $\add T_1\to\add T_2$. By Morita theory, there exists a progenerator $P\in\mod\Lambda$ such that $\fff$ is isomorphic to $\hom_{\Lambda}(P,\ )$ and $\endm_{\Lambda}(P)\simeq\Lambda$. Since $\Lambda$ is basic, $P$ is isomorphic to $\Lambda$ as a $\Lambda$-module. Thus we have an automorphism $\phi:\Lambda=\endm_{\Lambda}(\Lambda)\to\endm_{\Lambda}(P)\simeq\Lambda$. It is easily checked that $\fff$ is isomorphic to $\phi$. Since $M_i$ and $T_i$ are basic, we obtain the assertion.\rule{5pt}{10pt}

\vskip.5em{\bf\XE\ Non-commutative crepant resolution and representation dimension}

Throughout this section, fix a complete regular local ring $R$ of dimension $d\ge0$, an $R$-order $\Lambda$ which is an isolated singularity, and an $m$-cotilting $\Lambda$-module $T$. Put $\aa:=\mod\Lambda$ and $\bb:={}^\perp T$ as in \XCBB.

\vskip.5em{\bf\XEA\ }Let us start with studying properties of a $\Lambda$-module $M$ in terms of $\endm_\Lambda(M)$.

\vskip.5em{\bf Theorem }{\it
Let $M\in\bb$, $\Gamma:=\endm_\Lambda(M)$ and $n\ge 1$. Assume $\Lambda\oplus T\in\add M$.

(1) Assume $n<d$. Then $M$ is $(n-1)$-orthogonal if and only if $\depth_R\Gamma\ge n+1$.

(2) Assume that $M$ is $(m-1)$-orthogonal. Then $M$ is $(n-1)$-orthogonal if and only if $\Gamma$ satisfies the two-sided $(m+1,n+1)$-condition.

(3) Assume that $M$ is $(n-1)$-orthogonal. If $M\in\bb$ is maximal $(n-1)$-orthogonal, then $\gl\Gamma\le\max\{m,n+1\}$ holds, and the converse holds if $m\le n+1$.}

\vskip.5em{\sc Proof }(1) follows by \XCAA. By \XDD, $\dn{B}_{m,n}\subseteq\{(\Gamma,P,I)\in\dn{B}_{m,n}^q\ |\ \gl\Gamma\le\max\{n+1,m\}\}$ holds, and the equality holds if $m\le n+1$. Thus (3) follows by \XDDA. We will show (2). The `only if' part follows by \XDCA(5). To show the `if' part, we can assume $m\le n$. For $(\Lambda,M,T)\in\dn{A}_{m,m}^q$, put $(\Gamma,P,I):=\alpha(\Lambda,M,T)\in\dn{B}_{m,m}^q$. Since $\Gamma$ is an $R$-order and satisfies the two-sided $(m+1,n+1)$-condition, \XDE\ implies $(\Gamma,P,I)\in\dn{B}_{m,n}^q$. Thus $(\Lambda,M,T)\in\dn{A}_{m,n}^q$ holds, and we obtain $M\perp_{n-1}M$.\rule{5pt}{10pt}

\vskip.5em{\bf\XEB\ Definition }
Let us generalize the concept of Van den Bergh's non-commutative crepant resolution [V1,2] of commutative normal Gorenstein domains to our situation.

Again let $\Lambda$ be an $R$-order which is an isolated singularity. We say that $M\in\cm\Lambda$ gives a {\it Cohen-Macaulay non-commutative crepant resolution} ({\it CM NCCR} for short) $\Gamma:=\endm_\Lambda(M)$ of $\Lambda$ if $\Lambda\oplus D_d\Lambda\in\add M$ and $\Gamma$ is an $R$-order with $\gl\Gamma=d$. 

Our definition is slightly stronger than original non-commutative crepant resolutions in [V2] where $M$ is assumed to be reflexive (not Cohen-Macaulay) and $\Lambda\oplus D_d\Lambda\in\add M$ is not assumed. But all examples of non-commutative crepant resolutions in [V1,2] satisfy our condition. For the case $d\ge2$, we have the remarkable relationship below between CM NCCR and maximal $(d-2)$-orthogonal subcategories. In this case, $\Gamma$ is an Auslander algebra of type $(d,d,d-1)$ and has remarkable properties (see \XCBC\ and \XDG).

\vskip.5em{\bf\XEBA\ Theorem }{\it
Let $d\ge2$. Then $M\in\cm\Lambda$ gives a CM NCCR of $\Lambda$ if and only if $M\in\cm\Lambda$ is maximal $(d-2)$-orthogonal.}

\vskip.5em{\sc Proof }
$\Lambda\oplus D_d\Lambda\in\add M$ holds. Put $m:=d$ and $n:=d-1$ in \XEA(1) and (3).\rule{5pt}{10pt}

\vskip.5em{\bf\XEBB\ Example }
Let $k$ be a field of characteristic zero, $G$ a finite subgroup of $\Gl_d(k)$ with $d\ge2$, $\Omega:=k[[x_1,\cdots,x_d]]$ and $\Lambda:=\Omega^G$ the invariant subring. Assume that $G$ does not contain any pseudo-reflection except the identity, and that $\Lambda$ is an isolated singularity. In [I7;2.5], it is shown that $\cc:=\add{}_\Lambda\Omega$ is a maximal $(d-2)$-orthogonal subcategory of $\cm\Lambda$, and $\endm_\Lambda(\Omega)$ is the skew group ring $\Omega*G$ [A4] (see also [Y;10.8]). Hence $\Omega$ gives a CM NCCR $\Omega*G$ of $\Lambda$ (see [V2;1.1]), and $\Omega*G$ is an Auslander algebra of type $(d,d,d-1)$. We will study this example in \XFA.

\vskip.5em{\bf\XEC\ Conjecture }Fix a pair $(\aa,\bb)$ in \XCBB\ again, and $l\ge1$. It is interesting to study the relationship among all maximal $l$-orthogonal objects in $\bb$. Especially, we conjecture that {\it their endomorphism rings are derived equivalent}. It is suggestive to relate this conjecture to Van den Bergh's generalization [V2] of the Bondal-Orlov conjecture [BO], which asserts that {\it all (commutative or non-commutative) crepant resolutions of a normal Gorenstein domain have the same derived category}. Since maximal $l$-orthogonal subcategories are analogs of modules giving non-commutative crepant resolutions from the viewpoint of \XEBA, our conjecture is analogous to the Bondal-Orlov-Van den Bergh conjecture. We will show in \XECC\ that it is true for the case $l=1$.

\vskip.5em{\bf\XECA\ Lemma }{\it
Let $M_i\in\bb$ and $t\ge1$. Assume $\resdim{(\add M_1)}{M_2}\le t$, $M_2\perp_tM_2$, $M_1\perp_{t-1}M_2$ and that $M_1$ is a generator. Put $\Gamma_i:=\endm_\Lambda(M_i)$ and $U:=\hom_\Lambda(M_1,M_2)$. Then $U$ satisfies $\pd{}_{\Gamma_1}U\le t$, $({}_{\Gamma_1}U)\perp({}_{\Gamma_1}U)$ and $\endm_{\Gamma_1}(U)=\Gamma_2$.}

\vskip.5em{\sc Proof }Take a right $(\add M_1)$-resolution ${\bf X}:0\to X_t\to\cdots\to X_0\to M_2\to0$ of $M_2$, which is exact since $M_1$ is a generator. Then ${}_\Lambda(M_1,{\bf X}):0\to{}_\Lambda(M_1,X_t)\to\cdots\to{}_\Lambda(M_1,X_0)\to U\to0$ gives a projective resolution of the $\Gamma_1$-module $U$. Thus $\pd{}_{\Gamma_1}U\le t$ holds. By $M_2\perp_tM_2$ and $M_1\perp_{t-1}M_2$, we have an exact sequence ${}_\Lambda({\bf X},M_2):0\to\endm_\Lambda(M_2)\to{}_\Lambda(X_0,M_2)\to\cdots\to{}_\Lambda(X_t,M_2)\to0$. Since ${}_{\Gamma_1}({}_\Lambda(M_1,\ ),{}_\Lambda(M_1,M_2))={}_\Lambda(\ ,M_2)$ holds on $\add M_1$, the complex ${}_\Gamma({}_\Lambda(M_1,{\bf X}),U)$ is isomorphic to the exact sequence ${}_\Lambda({\bf X},M_2)$. Thus $({}_{\Gamma_1}U)\perp({}_{\Gamma_1}U)$ and $\endm_{\Gamma_1}(U)=\Gamma_2$ hold.\rule{5pt}{10pt}

\vskip.5em{\bf\XECB\ Theorem }{\it
Let $M_i\in\bb$ be maximal $l$-orthogonal ($i=1,2$) and $k\le l\le 2k+1$. Assume $M_1\perp_kM_2$. Put $\Gamma_i:=\endm_\Lambda(M_i)$ and $U:=\hom_\Lambda(M_1,M_2)$. Then $U$ is a tilting $(\Gamma_1,\Gamma_2)$-module with $\pd{}_{\Gamma_1}U\le l-k$ [M]. Thus $\Gamma_1$ and $\Gamma_2$ are derived equivalent.}

\vskip.5em{\sc Proof }
Put $t:=l-k$. By $M_1\perp_{k}M_2$ and \XBDA, $\resdim{(\add M_1)}{M_2}\le t$ and $\resdim{(\add M_2)^{op}}{M_1}\le t$ hold. We have $M_i\perp_{t}M_i$ and $M_1\perp_{t-1}M_2$. By \XECA, $\pd{}_{\Gamma_1}U\le t$, $({}_{\Gamma_1}U)\perp({}_{\Gamma_1}U)$ and $\endm_{\Gamma_1}(U)=\Gamma_2$ hold. Take a left $(\add M_2)$-resolution ${\bf Y}:0\to M_1\to Y_0\to\cdots\to Y_{t}\to0$ of $M_1$. By $M_1\perp_{t}M_1$ and $M_1\perp_{t-1}M_2$, ${}_\Lambda(M_1,{\bf Y}):0\to\Gamma_1\to{}_\Lambda(M_1,Y_0)\to\cdots\to{}_\Lambda(M_1,Y_{t})\to0$ is exact with ${}_\Lambda(M_1,Y_i)\in\add{}_{\Gamma_1}U$. Thus $U$ is a tilting $(\Gamma_1,\Gamma_2)$-module, and $\Gamma_1$ and $\Gamma_2$ are derived equivalent by a result of Happel [Ha].\rule{5pt}{10pt}

\vskip.5em{\bf\XECC\ Corollary }{\it
(1) Let $\cc_i=\add M_i$ be a maximal $1$-orthogonal subcategory of $\bb$ and $\Gamma_i:=\endm_\Lambda(M_i)$ ($i=1,2$). Then $\Gamma_1$ and $\Gamma_2$ are derived equivalent. In particular, $\#\ind\cc_1=\#\ind\cc_2$ holds.

(2) If $d\le 3$, then all CM NCCR of $\Lambda$ have the same derived category.}

\vskip.5em{\sc Proof }(1) Put $l:=1$ and $k:=0$ in \XECB. Since a derived equivalence preserves Grothendieck groups [Ha], the latter assertion follows.

(2) If $d=3$, then the assertion follows by (1) and \XEBA. If $d=2$, then any $M$ giving a CM NCCR satisfies $\cm\Lambda=\add M$ by \XEBA. Thus $\endm_\Lambda(M)$ is unique up to Morita-equivalences. Assume $d\le 1$. It is well-known that, if an $R$-order $\Gamma$ satisfies $\gl\Gamma=d$, then $\gl\endm_\Gamma(P)=d$ holds for any $P\in\add{}_\Gamma\Gamma$ [CR]. Thus $\Lambda$ has a CM NCCR if and only if $\gl\Lambda=d$. In this case, any CM NCCR is Morita-equivalent to $\Lambda$.\rule{5pt}{10pt}

\vskip.5em{\bf\XECD\ }We obtain the corollary below by \XBFB, \XECA\ and \XECB. For the case $m=0$, $\Gamma$ satisfies $\gl\Gamma\le 3$ and $\domdim\Gamma\ge3$ by \XDBA. Then Miyachi's theorem [Mi] implies that $\Omega^2(\mod\Gamma)$ coincides with the category of $\Gamma$-modules $X$ with $\pd{}_\Gamma X\le1$. Thus our corollary gives the (independent) result of Geiss-Leclerc-Schro\"er [GLS2].

\vskip.5em{\bf Corollary }{\it
Let $M\in\bb$ be maximal $1$-orthogonal and $\Gamma:=\endm_\Lambda(M)$. Assume $m\le 2$. Then we have equivalences $\fff:=\bb(M,\ ):\bb\to\Omega^2(\mod\Gamma)$ and $\ggg:=\bb(\ ,M):\bb\to\Omega^2(\mod\Gamma^{op})$, which send maximal $1$-orthogonal (resp. $1$-orthogonal) objects in $\bb$ to tilting (resp. partial tilting) $\Gamma$ and $\Gamma^{op}$-modules respectively and satisfy $\fff=(\ )^*\circ\ggg$ and $\ggg=(\ )^*\circ\fff$ for $(\ )^*=\hom_\Gamma(\ ,\Gamma)$.}

\vskip.5em{\bf\XED\ Definition }
Let us generalize the concept of Auslander's representation dimension [A1] to relate it to non-commutative crepant resolutios. For $n\ge1$, define the {\it $n$-th representation dimension} $\rdim_n\Lambda$ of an $R$-order $\Lambda$ which is an isolated singularity by
\[\rdim_n\Lambda:=\inf\{\gl\endm_\Lambda(M)\ |\ M\in\cm\Lambda,\ \Lambda\oplus D_d\Lambda\in\add M,\ M\perp_{n-1}M\}.\]
In other words, we consider all $(\Lambda,M,D_d\Lambda)\in\dn{A}_{d,n}^q$, and $\rdim_n\Lambda$ is the infimum of global dimension of corresponding quasi-Auslander algebras $\endm_\Lambda(M)$ of type $(d,d,n)$. Obviously $d\le\rdim_n\Lambda\le\rdim_{n^\prime}\Lambda$ holds for any $n\le n^\prime$ ([R;3.2]). Notice that $\rdim_1\Lambda$ coincides with the representation dimension defined in [A1][I6].

\vskip.5em{\bf\XEDA\ }We call $\Lambda$ {\it representation-finite} if $\#\ind(\cm\Lambda)<\infty$. In the sense of (1) below, $\rdim_1\Lambda$ measures how far $\Lambda$ is from being representation-finite (cf. [A1][I6]).

\vskip.5em{\bf Theorem }{\it 
(1) If $\Lambda$ is representation-finite, then $\rdim_1\Lambda\le\max\{2,d\}$ holds, and the converse holds if $d\le 2$. If $d>2$, then the converse does not necessarily hold.

(2) $\Lambda$ has a CM NCCR if and only if $\rdim_l\Lambda=d$ holds for $l:=\max\{1,d-1\}$.

(3) Let $n\ge1$. If $\cm\Lambda$ has a maximal $(n-1)$-orthogonal subcategory $\cc$ with $\#\ind\cc<\infty$, then $\rdim_n\Lambda\le\max\{n+1,d\}$ holds, and the converse holds if $d\le n+1$.}

\vskip.5em{\sc Proof }
(3) The assertion follows immediately by \XEA(3). 

(2) If $d\ge2$, then the assertion follows by (3) and \XEBA. For the case $d<2$, $\Lambda$ has a CM NCCR if and only if $\gl\Lambda=d$ if and only if $\rdim_1\Lambda=d$ by the argument in the proof of \XECC(2).

(1) We obtain the assertion by putting $n:=1$ in (3). Let us give a counter-example for $d>2$. Take $\Lambda$ in \XEBB, which has a CM NCCR. Thus $\rdim_1\Lambda=\rdim_{d-1}\Lambda=d$ by (2). But $\Lambda$ is representation-infinite except the case $d=3$ and $G=\langle{\rm diag}(-1,-1,-1)\rangle$ [AR2].\rule{5pt}{10pt}

\vskip.5em{\bf\XEDB\ }It is an interesting problem raised by Auslander [A1] to calculate the value of $\rdim_n\Lambda$. In particular, when $\rdim_n\Lambda$ is finite? For the case $n=1$ and $d\le1$, we have the finiteness result below obtained by the author [I1,3,6] recently (see also [L]).

\vskip.5em{\bf Theorem }{\it If $d\le 1$, then $\rdim_1\Lambda<\infty$.}

\vskip.5em{\bf\XEDC\ }
If $d\ge2$, then $\rdim_1\Lambda<\infty$ does not necessarily hold. For example, if $d=2$ and $\Lambda$ is representation-infinite commutative Gorenstein, then $\rdim_1\Lambda=\infty$ holds by (3) below, which we will prove by the argument of Van den Bergh in [V2;4.2]. We call a module-finite $R$-algebra $\Lambda$ {\it symmetric} if $D_d\Lambda$ is isomorphic to $\Lambda$ as a $(\Lambda,\Lambda)$-module. Any symmetric order is Gorenstein, and the converse holds if it is commutative.

\vskip.5em{\bf Proposition }{\it
Assume $d\ge2$ and that $\Lambda$ is a symmetric $R$-order.

(1) If $M\in\mod\Lambda$ is a reflexive $R$-module, then $\endm_\Lambda(M)$ is a symmetric $R$-algebra.

(2) $\rdim_{d-1}\Lambda$ is either $d$ or $\infty$.

(3) If $d=2$, then $\rdim_1\Lambda<\infty$ if and only if $\Lambda$ is representation-finite.}

\vskip.5em{\sc Proof }
(1) Put $\Gamma:=\endm_\Lambda(M)$. Since $\Lambda$ is symmetric, $M^*:={}_\Lambda(M,\Lambda)$ is isomorphic to $D_dM$ as a $(\Gamma,\Lambda)$-module. We have a natural map $f:M^*\otimes_\Lambda M\to\Gamma$. Since $\Lambda$ is assumed to be an isolated singularity, $f_{\dn{p}}$ is an isomorphism for any $\dn{p}\in\Spec R$ with $\height\dn{p}=1$. Now consider a $(\Gamma,\Gamma)$-homomorphism $D_d\Gamma\stackrel{D_df}{\to}D_d(M^*\otimes_\Lambda M)\simeq D_d(D_dM\otimes_\Lambda M)={}_\Lambda(M,D_dD_dM)=\Gamma$. Since $D_df$ is a map between reflexive $R$-modules such that $(D_df)_{\dn{p}}$ is an isomorphism for any $\dn{p}\in\Spec R$ with $\height\dn{p}=1$, it is an isomorphism.

(2) Take $M\in\cm\Lambda$ with $M\perp_{d-2}M$. Then $\Gamma:=\endm_\Lambda(M)$ is an $R$-order by \XEA(1). Since $\Gamma$ is Gorenstein by (1), $\id{}_\Gamma\Gamma=d$ holds by \XCAB. Thus $\gl\Gamma=d$ or $\infty$.

(3) The assertion follows by \XEDA(1).\rule{5pt}{10pt}

\vskip.5em{\bf\XEDD\ }
We end this subsection by giving a few remarks on the value of $\rdim_n\Lambda$.

(1) Assume $d=0$. Thus $\rdim_1\Lambda<\infty$ holds by \XEDB, and $\rdim_1\Lambda\le2$ if and only if $\Lambda$ is representation-finite by \XEDA. Dugas [D] and Guo [G] independently proved that $\rdim_1\Lambda$ is preserved by stable equivalences. Recently, many results are obtained on algebras with $\rdim_1\Lambda\le3$, for example, they satisfy the famous finitistic dimension conjecture [IT]. Many classes of algebras are known to satisfy $\rdim_1\Lambda\le3$, e.g. hereditary algebras [A1], tilted algebras [APT], algebras with radical square zero [A1], special biserial algebras [EHIS] and so on. See also [BHS][CP][Ho][X].

On the other hand, Rouquier showed that $\rdim_1\Lambda=l+1$ holds for the exterior algebra $\Lambda=\wedge(k^l)$ of the $l$-dimensional vector space by applying his concept of the dimension of triangulated categories [Ro1,2] (see also [KK]). In general, it seems to be difficult to know the precise value of $\rdim\Lambda$ when this is larger than $3$.

(2) Assume $d=2$. Then $\rdim_1\Lambda=2$ if and only if $\Lambda$ has a CM NCCR if and only if $\Lambda$ is representation-finite by \XEDA. If $\Lambda$ is commutative and contains its residue field $\ccc$, then it is equivalent to be a quotient singularity [A4].

(3) A trivial example of an order $\Lambda$ with $\rdim_n\Lambda=\infty$ is an order which does not satisfy $D_d\Lambda\perp_{n-1}\Lambda$. Let us give a non-trivial Gorenstein example. Van den Bergh proved that $\Lambda:=k[[x,y,z,t]]/(x^2+y^2+z^2+t^{2b+1})$ does not have a non-commutative crepant resolution [V1;A.1]. More strongly, he proved that there is no non-free reflexive $\Lambda$-module $M$ such that $\endm_\Lambda(M)$ is an order, or equivalently, $M\perp_1M$ holds (\XEA(1)). Since $\Lambda$ is not regular, $\rdim_2\Lambda=\infty$ holds.

\vskip.5em{\bf\XEE\ Conjecture }For $l\ge1$ and $\bb$ in \XCBB, it seems that no example of a maximal $l$-orthogonal subcategory $\cc$ of $\bb$ with $\#\ind\cc=\infty$ is known. This suggests us to study
\[ o(\bb):=\sup_{\cc\subseteq\bb,\ \cc\perp_1\cc}\#\ind\cc.\]
We conjecture that {\it $o(\bb)$ is always finite}. If $\Lambda$ is a preprojective algebra of Dynkin type $\Delta$, then Geiss-Schr\"oer [GS] proved that $o(\mod\Lambda)$ equals the number of positive roots of $\Delta$ (see \XFBA). It would be interesting to find an interpretation of $o(\bb)$ for more general $\bb$.

\vskip.5em{\bf\XEEA\ }
For some classes of $\bb$, one can calculate $o(\bb)$ by using the theorem below. Especially, (1) seems to be interesting in connection with known results in \XEDD(1).

\vskip.5em{\bf Theorem }{\it
(1) $\rdim_1\Lambda\le3$ implies $o(\cm\Lambda)<\infty$.

(2) If $\bb$ has a maximal $1$-orthogonal subcategory $\cc$, then $o(\bb)=\#\ind\cc$.

(3) If $\bb$ has a subcategory $\cc$ such that $\Lambda\in\cc$ and $\resdim{\cc}{\bb}\le 1$, then $o(\bb)\le\#\ind\cc$.}

\vskip.5em{\sc Proof }
We will show (3). We can assume $\cc=\add{}_\Lambda M_1$. For any $1$-orthogonal $M_2\in\bb$ and $t:=1$, we apply \XECA. Consequently, $U$ is a partial tilting module. Since any partial tilting module is a direct summand of a tilting module [Ha], we obtain $\#\ind(\add{}_\Lambda M_2)=\#\ind(\add{}_{\Gamma_1}U)\le\#\ind(\add{}_{\Gamma_1}\Gamma_1)=\#\ind\cc$. Thus $o(\bb)\le\#\ind\cc$ holds. In particular, (2) follows. We will show (1). Take $M\in\cm\Lambda$ such that $\Lambda\oplus D_d\Lambda\in\add M$ and $\gl\endm_\Lambda(M)\le3$. It is easily shown that $\resdim{(\add M)}{(\cm\Lambda)}\le1$ holds (e.g. [EHIS;2.1]). By (3), $o(\cm\Lambda)\le\#\ind(\add M)$ holds.\rule{5pt}{10pt}

\vskip.5em{\bf\XEEB\ }Concerning our conjecture, let us recall the well-known proposition below which follows by a geometric argument due to Voigt ([P;4.2]). It is interesting to ask whether it is true without the restriction on $R$. If it is true, then any $1$-orthogonal subcategory of $\bb$ is `discrete', and our conjecture asserts that it is finite. It is interesting to study the discrete structure of $1$-orthogonal objects in $\bb$ and the relationship to the whole structure of $\bb$.

\vskip.5em{\bf Proposition }{\it
Assume $d=0$ and that $R$ is an algebraically closed field. For any $n>0$, there are only finitely many isoclasses of $1$-orthogonal $\Lambda$-modules $X$ with $\dim_RX=n$.}

\vskip.5em{\bf\XF\ Applications and examples }

\vskip.5em{\bf\XFA\ }
Let us recall Auslander's contribution to McKay correspondence [Mc]. He proved in [A4] (see also [Y]) that the McKay quiver of a finite subgroup $G$ of $\Gl_2(k)$ coincides with the Auslander-Reiten quiver of the invariant subring $k[[x,y]]^G$. The aim of this section is to give a higher dimensional generalization \XFAD\ of this result.

\vskip.5em{\bf\XFAA\ Definiton }
Let $(\aa,\bb)$ be a pair in \XCBB\ and $\cc$ a maximal $(n-1)$-orthogonal subcategory of $\bb$. We will define the {\it Auslander-Reiten quiver} $\dn{A}(\cc)$ of $\cc$. For simplicity, we assume that the residue field $k$ of $R$ is algebraically closed. The set of vertices of $\dn{A}(\cc)$ is $\ind\cc$. For $X,Y\in\ind\cc$, we denote by $d_{XY}$ be the multiplicity of $X$ in $C$ for the sink map $C\to Y$ (\XAA), which equals to the multiplicity of $Y$ in $C^\prime$ for the source map $X\to C^\prime$. Draw $d_{XY}$ arrows from $X$ to $Y$. Draw a dotted arrow from non-projective $X\in\ind\cc$ to non-injective $\tau_nX\in\ind\cc$. If $\cc=\add M$, then $\dn{A}(\cc)$ coincides with the Gabriel quiver of $\endm_\Lambda(M)$ since $d_{XY}=\dim_kJ_{\cc}/J_{\cc}^2(X,Y)$ holds.

\vskip.5em{\bf\XFAB\ Definition }Let $k$ be a field of characteristic zero and $G$ a finite subgroup of $\Gl_d(k)$ with $d\ge2$. Recall that the {\it McKay quiver} $\dn{M}(G)$ of $G$ [Mc] is defined as follows: The set of vertices is the set $\irr G$ of isoclasses of irreducible representations of $G$. Let $V$ be the representation of $G$ acting on $k^d$ through $\Gl_d(k)$. For $X,Y\in\irr G$, we denote by $d_{XY}$ the multiplicity of $X$ in $V\otimes_kY$, and draw $d_{XY}$ arrows from $X$ to $Y$. Let $S=\wedge^dV$ be the $1$-dimensional representation of $G$ given by the determinant. Draw a dotted arrow from $X\in\irr G$ to $\tau_{d-1}X:=S\otimes_kX\in\irr G$.

\vskip.5em{\bf\XFAC\ }
Let $G$ be in \XFAB, $\Omega:=k[[x_1,\cdots,x_d]]$, $\Lambda:=\Omega^G$ the invariant subring and $\Gamma:=\Omega*G$ the skew group ring. Assume that $G$ does not contain any pseudo-reflection except the identity, and that $\Lambda$ is an isolated singularity. Then $\cc:=\add{}_\Lambda\Omega$ forms a maximal $(d-2)$-orthogonal subcategory of $\cm\Lambda$ and $\endm_\Lambda(\Omega)=\Gamma$ holds by \XEBB. Let us compare $\dn{A}(\cc)$ and $\dn{M}(G)$ by applying the argument in [A4] (see also [Y]) to arbitrary $d\ge2$.

The functor $\fff:=\Omega\otimes_k\ :\mod kG\to\mod\Gamma$ induces a bijection $\irr G\to\ind(\add{}_\Gamma\Gamma)$ [A4] (see also [Y;10.1]). Define a functor $\ggg:\mod\Gamma\to\mod\Lambda$ by $\ggg(X):=X^G$ and $\ggg(f):=f|_{X^G}$ for $X,Y\in\mod\Gamma$ and $f\in\hom_\Gamma(X,Y)$. Since $\endm_\Lambda(\Omega)=\Gamma$ holds, $\ggg$ restricts to the equivalence $\ggg:\add{}_{\Gamma}\Gamma\to\add{}_\Lambda\Omega=\cc$. Composing $\fff$ and $\ggg$, we obtain a functor $\hhh:=\ggg\circ\fff:\mod kG\to\cc$, which gives a bijection $\hhh:\irr G\to\ind\cc$. Let
\[{\bf K}:0\to\Omega\otimes_k\wedge^dV\to\cdots\to\Omega\otimes_k\wedge^2V\to\Omega\otimes_k V\to \Omega\to k\to0\]
be the Koszul complex of $\Omega$. Then ${\bf K}$ forms an exact sequence of $\Gamma$-modules. For any $X\in\irr G$,
\[{\bf K}\otimes_kX:0\to\fff(\wedge^dV\otimes_kX)\to\cdots\to\fff(\wedge^2V\otimes_kX)\to\fff(V\otimes_kX)\to\fff(X)\to X\to0\]
gives a minimal projective resolution of the $\Gamma$-modules $X$. Taking $\ggg$, we obtain an exact sequence
\[\ggg({\bf K}\otimes_kX):0\to\hhh(\wedge^dV\otimes_kX)\stackrel{f_{d-1}}{\to}\cdots\stackrel{f_2}{\to}\hhh(\wedge^2V\otimes_kX)\stackrel{f_1}{\to}\hhh(V\otimes_kX)\stackrel{f_0}{\to}\hhh(X)\to\ggg(X)\to0\]
with $f_i\in J_{\cc}$ for any $i$. Now $\ggg(X)=X$ holds if $X$ is a trivial $G$-module, and $\ggg(X)=0$ otherwise. Since $\ggg:\add{}_{\Gamma}\Gamma\to\cc$ was an equivalence, $0\to\cc(\ ,\hhh(\wedge^dV\otimes_kX))\stackrel{\cdot f_{d-1}}{\to}\cdots\stackrel{\cdot f_1}{\to}\cc(\ ,\hhh(V\otimes_kX))\stackrel{\cdot f_0}{\to}J_{\cc}(\ ,\hhh(X))\to0$ is exact. Consequently, $\ggg({\bf K}\otimes_kX)$ is a $(d-1)$-fundamental sequence (\XCBC) if $X$ is trivial, and a $(d-1)$-almost split sequence (\XBEC) otherwise. Thus we obtain the theorem below, where we put $\tau_{d-1}\Lambda:=\nu_\Lambda\Lambda=D_d\Lambda$ from the viewpoint of \XCBC.

\vskip.5em{\bf\XFAD\ Theorem }{\it
The Auslander-Reiten quiver $\dn{A}(\cc)$ of $\cc:=\add{}_\Lambda\Omega$ coincides with the McKay quiver $\dn{M}(G)$ of $G$. Precisely speaking, there exists a bijection $\hhh:\irr G\to\ind\cc$ such that $\hhh\circ\tau_{d-1}=\tau_{d-1}\circ\hhh$ and $d_{XY}=d_{\hhh(X),\hhh(Y)}$ for any $X,Y\in\irr G$.}

\vskip.5em{\bf\XFB\ }
Geiss-Leclerc-Schr\"oer [GLS1,2] applied $1$-orthogonal (={\it rigid} in their papers) modules to study semicanonical basis of the quantum enveloping algebra. Their work is closely related to our study in this paper. Let $\Lambda$ be a preprojective algebra of type $A_n$ over an algebraically closed field $k$. Thus $\Lambda$ is defined by the following quiver with relations $a_1b_1=0$, $a_{i+1}b_{i+1}=b_ia_i$ ($1\le i\le n-2$) and $b_{n-1}a_{n-1}=0$.
\[\stackrel{e_1}{\bullet}\def\arraystretch{.5}\begin{array}{c}{\scriptscriptstyle a_1}\\ \longrightarrow\\ \longleftarrow\\ {\scriptscriptstyle b_1}\end{array}\stackrel{e_2}{\bullet}\begin{array}{c}{\scriptscriptstyle a_2}\\ \longrightarrow\\ \longleftarrow\\ {\scriptscriptstyle b_2}\end{array}\stackrel{e_3}{\bullet}\begin{array}{c}{\scriptscriptstyle a_3}\\ \longrightarrow\\ \longleftarrow\\ {\scriptscriptstyle b_3}\end{array}\cdots\begin{array}{c}{\scriptscriptstyle a_{n-1}}\\ \longrightarrow\\ \longleftarrow\\ {\scriptscriptstyle b_{n-1}}\end{array}\stackrel{e_n}{\bullet}\]
We denote by $\iota$ the automorphism of $\Lambda$ defined by $\iota(e_i):=e_{n+1-i}$, $\iota(a_i)=b_{n-i}$ and $\iota(b_i)=a_{n-i}$ for any $i$. Let us collect results on $1$-orthogonal subcategories of $\mod\Lambda$. Especially, (4) below answers a question raised by Schr\"oer.

\vskip.5em{\bf\XFBA\ Theorem }{\it
Let $\Lambda$ be a preprojective algebra of type $A_n$ over an algebraically closed field $k$.

(1)(Geiss-Leclerc-Schr\"oer) Any $1$-orthogonal subcategory of $\mod\Lambda$ is contained in a maximal $1$-orthogonal subcategory of $\mod\Lambda$.

(2) Any maximal $1$-orthogonal subcategory $\cc$ of $\mod\Lambda$ satisfies $\#\ind\cc=\frac{n(n+1)}{2}$.

(3) Let $M\in\mod\Lambda$ be a generator and $\Gamma:=\endm_\Lambda(M)$. Then $M$ is (maximal) $1$-orthogonal if and only if $\domdim\Gamma\ge 3$ (and $\gl\Gamma\le 3$).

(4) Let $M$ and $M^\prime\in\mod\Lambda$ be basic $1$-orthogonal generators. Then $\endm_\Lambda(M)$ is isomorphic to $\endm_\Lambda(M^\prime)$ if and only if $M$ is isomorphic to $M^\prime$ or $\iota(M^\prime)$.}

\vskip.5em{\sc Proof }
Geiss-Leclerc-Schr\"oer proved (1) in [GLS2]. They constructed maximal $1$-orthogonal $M\in\mod\Lambda$ with $\#\ind(\add M)=\frac{n(n+1)}{2}$. Thus (2) follows by \XECC(1). (3) follows by \XEA(2)(3). We will show (4) in the rest of this section.\rule{5pt}{10pt}

\vskip.5em{\bf\XFBB\ }Let us start with calculating the group $\aut(\Lambda)$ of $k$-algebra automorphisms.

\vskip.5em{\bf Proposition }{\it
Put $H:=\{g\in\aut(\Lambda)\ |\ g$ fixes any $e_i$ and $a_i\}$. Then $\aut(\Lambda)=(\inn(\Lambda)\timesl\langle\iota\rangle)\timesl H$ and $H\simeq\aut(k[x]/(x^{l+1}))$ hold, where $l$ is the maximal integer which does not exceed $n/2$.}

\vskip.5em{\sc Proof }
(i) Let $\{f_1,\cdots,f_n\}$ be a complete set of orthogonal primitive idempotents of $\Lambda$. We will show that there exist $\lambda\in\Lambda^\times$ and $\sigma\in\dn{S}_n$ such that $e_i\lambda=\lambda f_{\sigma(i)}$ for any $i$.

Since $\bigoplus_{i=1}^n\Lambda e_i=\Lambda=\bigoplus_{i=1}^n\Lambda f_i$ holds, Krull-Schmidt theorem implies that there exists $\sigma\in\dn{S}_n$ such that $\Lambda e_i\simeq\Lambda f_{\sigma(i)}$ holds for any $i$. Since $\hom_\Lambda(\Lambda e_i,\Lambda f_{\sigma(i)})$ (resp. $\hom_\Lambda(\Lambda f_{\sigma(i)},\Lambda e_i)$) can be identified with $e_i\Lambda f_{\sigma(i)}$ (resp. $f_{\sigma(i)}\Lambda e_i$), there exist $\lambda_i\in e_i\Lambda f_{\sigma(i)}$ and $\gamma_i\in f_{\sigma(i)}\Lambda e_i$ such that $\lambda_i\gamma_i=e_i$ and $\gamma_i\lambda_i=f_{\sigma(i)}$. Put $\lambda:=\sum_{i=1}^n\lambda_i\in\Lambda$ and $\gamma:=\sum_{i=1}^n\gamma_i\in\Lambda$. Then $\lambda\gamma=1=\gamma\lambda$ holds by $\lambda\in\Lambda^\times$, and $e_i\lambda=\lambda_i=\lambda f_{\sigma(i)}$ holds.

(ii) We will show that $\aut(\Lambda)$ is generated by $\inn(\Lambda)$, $\iota$ and $H$.

Put $G:=\{g\in\aut(\Lambda)\ |\ g$ fixes any $e_i\}$ and fix $g\in\aut(\Lambda)$. By (i), there exist $h\in\inn(\Lambda)$ and $\sigma\in\dn{S}_n$ such that $hg(e_i)=e_{\sigma(i)}$ for any $i$. It is easily shown that $\sigma$ is either identity or $\sigma(i)=n+1-i$ for any $i$. Thus $hg\in G$ holds for the former case, and $\iota hg\in G$ holds for the latter case. Now we only have to show that $G$ is generated by $\inn(\Lambda)$ and $H$. Again fix $g\in G$. Inductively, we can take $\lambda_i\in(e_i\Lambda e_i)^\times$ such that $\lambda_1:=e_1$ and $\lambda_ig(a_i)=a_i\lambda_{i+1}$ for any $i$, since $e_i\Lambda e_{i+1}$ is generated by $a_i$ as an $(e_{i+1}\Lambda e_{i+1})^{op}$-module. Put $\lambda:=\sum_{i=1}^n\lambda_i\in\Lambda^\times$. Then $\lambda e_i\lambda^{-1}=e_i$ and $\lambda g(a_i)\lambda^{-1}=a_i$ hold for any $i$. Thus $(\lambda\cdot\lambda^{-1})g\in H$ holds.

(iii) We will show the latter equality. Fix $g\in\aut(k[x]/(x^{l+1}))$. Then there exist $s_1\in k^\times$ and $s_j\in k$ ($1<j\le l$) such that $g(x)=\sum_{j=1}^ls_jx^j$. Define $\phi(g)\in H$ by $\phi(g)(b_i):=\sum_{j=1}^ls_j(b_ia_i)^{j-1}b_i$ for any $i$. Thus we have a map $\phi:\aut(k[x]/(x^{l+1}))\to H$, which is easily checked to be an injective homomorphism. We will show that $\phi$ is surjective. Fix $h\in H$. By the relation $a_{i+1}h(b_{i+1})=h(b_i)a_i$, it is easily checked that there exist $s_1\in k^\times$ and $s_j\in k$ ($1<j\le l$) such that $h(b_i)=\sum_{j=1}^ls_j(b_ia_i)^{j-1}b_i$ for any $i$. Define $g\in\aut(k[x]/(x^{l+1}))$ by $g(x)=\sum_{j=1}^ls_jx^j$. Then $\phi(g)=h$ holds.

(iv) We will show the former equality. We can assume $n>1$. Take $g\in\inn(\Lambda)\cap H$ and put $g=(\lambda\cdot\lambda^{-1})$ for $\lambda\in\Lambda^\times$. Since $\lambda e_i=e_i\lambda$ holds for any $i$, we can put $\lambda=\sum_{i=1}^n\lambda_i$ with $\lambda_i\in(e_i\Lambda e_i)^\times$. Since $\lambda_ia_i=a_i\lambda_{i+1}$ holds for any $i$, we can easily check that $\lambda\in\cen\Lambda$ and $g=1$ hold. Take $g\in\inn(\Lambda)\iota\cap H$ and put $g=(\lambda\cdot\lambda^{-1})\iota$ for $\lambda\in\Lambda^\times$. Then $\lambda e_{n+1-i}=e_i\lambda$ holds for any $i$. Thus $e_1\lambda e_1=0$ holds by $n>1$, a contradiction to $e_1(\Lambda^\times)e_1=(e_1\Lambda e_1)^\times$. Consequently, $\inn(\Lambda)\cap\langle\iota\rangle=1$ and $(\inn(\Lambda)\timesl\langle\iota\rangle)\cap H=1$ hold. We only have to check that $\inn(\Lambda)\timesl\langle\iota\rangle$ is normalized by $H$. We will show that $h^{-1}\iota h\iota\in\inn(\Lambda)$ holds for any $h\in H$. The isomorphism $\phi$ in (iii) shows that there exists $c\in(\cen\Lambda)^\times$ such that $h(b_i)=b_ic$ for any $i$. Put $d:=h^{-1}(c)\in(\cen\Lambda)^\times$. Since $\iota$ acts trivially on $\cen\Lambda$ and $h^{-1}(b_i)=b_id^{-1}$ holds for any $i$, we have $h^{-1}\iota h\iota(a_i)=a_id$ and $h^{-1}\iota h\iota(b_i)=b_id^{-1}$ for any $i$. It is easily checked that $\lambda:=\sum_{i=1}^n e_id^i$ satisfies $h^{-1}\iota h\iota=(\lambda\cdot\lambda^{-1})\in\inn(\Lambda)$.\rule{5pt}{10pt}

\vskip.5em{\bf\XFBC\ Lemma }{\it
Let $\Lambda$ be a finite-dimensional $k$-algebra and $d>0$. Put $F:=\{g\in\aut(\Lambda)\ |\ g(X)$ is isomorphic to $X$ for any $1$-orthogonal $\Lambda$-module $X$ with $\dim_kX=d\}$. Then $F$ forms a Zariski-closed subgroup of $\aut(\Lambda)$ of finite index.}

\vskip.5em{\sc Proof }
Fix a $k$-basis $\{x_1\cdots,x_n\}$ of $\Lambda$. We denote by $\mod_d\Lambda$ the set of $\Lambda$-module structure on the $k$-vector space $k^d$. Denoting the action of $x_i$ on $k^d$ by $X_i\in\ma_d(k)\simeq{\bf A}_k^{d^2}$, we regard $\mod_d\Lambda$ as a closed subset of ${\bf A}_k^{d^2n}$. Thus $\Gl_d(k)$ acts on $\mod_d\Lambda$ by $g(X_1,\cdots,X_n):=(gX_1g^{-1},\cdots,gX_ng^{-1})$, and each $\Gl_d(k)$-orbit corresponds to an isoclass of $\Lambda$-modules. On the other hand, the closed subgroup $\aut(\Lambda)$ of $\Gl_n(k)$ acts on $\mod_d\Lambda$ by $a(X_1,\cdots,X_n):=(\sum_{j=1}^na_{1j}X_j,\cdots,\sum_{j=1}^na_{nj}X_j)$ for $a=(a_{ij})_{1\le i,j\le n}$.

Put $E:=\{X\in\mod_d\Lambda\ |\ \ext^1_\Lambda(X,X)=0\}$. By Voigt's lemma, there are only finitely many $\Gl_d(k)$-orbits $C_1,\cdots,C_t$ in $E$, and each $C_i$ forms a connected component of $E$. Since the action of $\Gl_d(k)$ and $\aut(\Lambda)$ on $\mod_d\Lambda$ commute, $\aut(\Lambda)$ acts on the set of orbits $S:=\{C_1,\cdots,C_t\}$. Since $F$ is the kernel of this action, our assertions follow.\rule{5pt}{10pt}

\vskip.5em{\bf\XFBD\ Proof of \XFBA(4) }
Fix a basic $1$-orthogonal generator $X$. By \XDH, we only have to show that $g(X)$ is isomorphic to $X$ or $\iota(X)$ for any $g\in\aut(\Lambda)$. Put $F^\prime:=\{g\in\aut(\Lambda)\ |\ g(X)$ is isomorphic to $X$ or $\iota(X)\}$. Then $F^\prime$ contains $\inn(\Lambda)\timesl\langle\iota\rangle$. By \XFBC, $F^\prime$ forms a closed subgroup of $\aut(\Lambda)$ of finite index. By \XFBB, $\aut(\Lambda)/(\inn(\Lambda)\timesl\langle\iota\rangle)=H\simeq\aut(k[x]/(x^{l+1}))=k^\times\times{\bf A}_k^l$ is connected. Hence $F^\prime=\aut(\Lambda)$ holds.\rule{5pt}{10pt}

\vskip.5em{\bf\XFC\ }
Let $\Lambda$ be a representation-finite selfinjective algebra. It is well-known that the stable Auslander-Reiten quiver of $\Lambda$ has the form $\dn{A}(\underline{\mod}\Lambda)=\zzz\Delta/G$ for a Dynkin diagram $\Delta$ and an automorphism group $G$ of $\zzz\Delta$ [Ri]. In [I7], maximal $1$-orthogonal subcategories of $\mod\Lambda$ are characterized combinatorially in terms of $\zzz\Delta$. Let us study them from a little bit different viewpoint.

Let $H$ be a hereditary $k$-algebra of Dynkin type $\Delta$, $\dd:=D^b(\mod H)$ the bounded derived category and $F:=\tau^{-1}\circ[1]$ the auto-equivalence of $\dd$. We put $\ext^1_{\dd}(X,Y):=\hom_{\dd}(X,Y[1])$ for $X,Y\in\dd$. Buan-Marsh-Reineke-Reiten-Todorov [BMRRT] introduced the {\it cluster category} $\cc_H:=\dd/F$ and showed that $\cc_H$ is closely related to cluster algebras of Fomin-Zelevinsky [FZ1,2]. A key concept was an {\it $\ext$-configuration}, which is a subset $\tt$ of $\ind\dd$ satisfying $\tt=\{X\in\ind\dd\ |\ \ext^1_{\dd}(\tt,X)=0\}$. Since this definition is essentially similar to our maximal $1$-orthogonal subcategories in \XBD, we can regard $\ext$-configurations as `maximal $1$-orthogonal subcategories of triangulated categories'. 


\vskip.5em{\bf\XFCA\ }It is well known [Ha] that the bounded derived category $D^b(\mod H)$ is equivalent to the mesh category $k(\zzz\Delta)$ [Ri]. We identify $D^b(\mod H)$ with $k(\zzz\Delta)$.

\vskip.5em{\bf Theorem }{\it
Let $\Delta$ be a Dynkin diagram and $\Lambda$ a standard selfinjective algebra with a covering functor $\ppp:k(\zzz\Delta)\to\underline{\mod}\Lambda$. 

(1) The following diagrams are commutative.
{\small\[\begin{diag}
D^b(\mod H)&\RA{[1]}&D^b(\mod H)&\ \ \ \ \ &D^b(\mod H)&\RA{F}&D^b(\mod H)\\ 
\downarrow^{\ppp}&&\downarrow^{\ppp}&&\downarrow^{\ppp}&&\downarrow^{\ppp}\\ 
\underline{\mod}\Lambda&\RA{\Omega^-}&\underline{\mod}\Lambda&&\underline{\mod}\Lambda&\RA{\tau_2^-}&\underline{\mod}\Lambda
\end{diag}\]}

\vskip-1em
(2) For a subcategory $\cc$ of $\mod\Lambda$ containing $\Lambda$, the conditions below are equivalent.

\strut\kern1em
(i) $\cc$ is a maximal $1$-orthogonal subcategory of $\mod\Lambda$.

\strut\kern1em
(ii) $\ppp^{-1}(\ind\underline{\cc})$ is an $\ext$-configuration.}

\vskip.5em{\sc Proof }
(1) Since $\ppp$ is a triangle functor and $\Omega^-:\underline{\mod}\Lambda\to\underline{\mod}\Lambda$ and $[1]:D^b(\mod H)\to D^b(\mod H)$ are shift functors, the left diagram is commutative. Since $\ppp$ commutes with $\tau$, we obtain $\tau_2^-\circ\ppp=\tau^{-1}\circ\Omega^-\circ\ppp=\ppp\circ\tau^{-1}\circ[1]=\ppp\circ F$.

(2) We can take an automorphism group $G$ of $\zzz\Delta$ such that $\dd/G=k(\zzz\Delta)/G\simeq\underline{\mod}\Lambda$. We have an isomorphism $\ext^1_\Lambda(\ppp X,\ppp Y)\simeq\coprod_{g\in G}\ext^1_{\dd}(gX,Y)$ for any $X,Y\in\dd$ [BMRRT]. Thus $\ppp X\perp_1\ppp Y$ if and only if $GX\perp_1Y$ if and only if $GX\perp_1GY$. This implies $\ppp^{-1}(\underline{\cc}^{\perp_1})=(\ppp^{-1}\underline{\cc})^{\perp_1}$. Thus the assertion follows.\rule{5pt}{10pt}

\vskip.5em{\footnotesize
\begin{center}
{\bf References}
\end{center}

[Ar] M. Artin: Maximal orders of global dimension and Krull dimension two. Invent. Math. 84 (1986), no. 1, 195--222.


[ArS] M. Artin, W. F. Schelter: Graded algebras of global dimension $3$. Adv. in Math. 66 (1987), no. 2, 171--216.

[APT] I. Assem, M. I. Platzeck, S. Trepode: On the representation dimension of tilted and laura algebras. J. Algebra 296 (2006), no. 2, 426--439.

[A1] M. Auslander: Representation dimension of Artin algebras. Lecture notes, Queen Mary College, London, 1971.

[A2] M. Auslander: Functors and morphisms determined by objects. Representation theory of algebras (Proc. Conf., Temple Univ., Philadelphia, Pa., 1976), pp. 1--244. Lecture Notes in Pure Appl. Math., Vol. 37, Dekker, New York, 1978.

[A3] M. Auslander: Isolated singularities and existence of almost split sequences. Representation theory, II (Ottawa, Ont., 1984), 194--242, Lecture Notes in Math., 1178, Springer, Berlin, 1986.

[A4] M. Auslander: Rational singularities and almost split sequences. Trans. Amer. Math. Soc. 293 (1986), no. 2, 511--531. 


[AB] M. Auslander, R. Buchweitz: The homological theory of maximal Cohen-Macaulay approximations. Colloque en l'honneur de Pierre Samuel (Orsay, 1987). Mem. Soc. Math. France (N.S.) No. 38, (1989), 5--37.

[AR1] M. Auslander, I. Reiten: Stable equivalence of dualizing $R$-varieties. Advances in Math. 12 (1974), 306--366. 

[AR2] M. Auslander, I. Reiten: The Cohen-Macaulay type of Cohen-Macaulay rings. Adv. in Math. 73 (1989), no. 1, 1--23.

[AR3] M. Auslander, I. Reiten: Applications of contravariantly finite subcategories. Adv. Math. 86 (1991), no. 1, 111--152.

[AR4] M. Auslander, I. Reiten: $k$-Gorenstein algebras and syzygy modules. J. Pure Appl. Algebra 92 (1994), no. 1, 1--27.


[AR5] M. Auslander, I. Reiten: $D$Tr-periodic modules and functors. Representation theory of algebras (Cocoyoc, 1994), 39--50, CMS Conf. Proc., 18, Amer. Math. Soc., Providence, RI, 1996.

[ARS] M. Auslander, I. Reiten, S. O. Smalo: Representation theory of Artin algebras. Cambridge Studies in Advanced Mathematics, 36. Cambridge University Press, Cambridge, 1995.

[ARo] M. Auslander, K. W. Roggenkamp: A characterization of orders of finite lattice type. Invent. Math. 17 (1972), 79--84.

[AS] M. Auslander, S. O. Smalo: Almost split sequences in subcategories. J. Algebra 69 (1981), no. 2, 426--454. 

[ASo] M. Auslander, O. Solberg: Gorenstein algebras and algebras with dominant dimension at least $2$. Comm. Algebra 21 (1993), no. 11, 3897--3934.

[B] J.-E. Bjork: The Auslander condition on Noetherian rings. Seminaire d'Algebre Paul Dubreil et Marie-Paul Malliavin, 39eme Annee (Paris, 1987/1988), 137--173, Lecture Notes in Math., 1404, Springer, Berlin, 1989. 

[BHS] R. Bocian, T. Holm, A. Skowro\'nski: The representation dimension of domestic weakly symmetric algebras. Cent. Eur. J. Math. 2 (2004), no. 1, 67--75.

[BMRRT] A. Buan, R. Marsh, M. Reineke, I. Reiten, G. Todorov: Tilting theory and cluster combinatorics, to appear in Adv. Math.

[BO] A. Bondal, D. Orlov: Semiorthogonal decomposition for algebraic varieties, preprint arXiv:alg-geom/9506012.

[CR] C. W. Curtis, I. Reiner: Methods of representation theory. Vol. I. With applications to finite groups and orders. Wiley Classics Library. A Wiley-Interscience Publication. John Wiley \& Sons, Inc., New York, 1990.

[C] J. Clark: Auslander-Gorenstein rings for beginners. International Symposium on Ring Theory (Kyongju, 1999), 95--115, 
Trends Math., Birkhauser Boston, Boston, MA, 2001. 

[CP] F. U. Coelho, M. I. Platzeck: On the representation dimension of some classes of algebras. J. Algebra 275 (2004), no. 2, 615--628.

[D] A. Dugas: Representation dimension as a relative homological invariant of stable equivalence, to appear in Algebr. Represent. Theory.

[EHIS] K. Erdmann, T. Holm, O. Iyama, J. Schr\"oer: Radical embeddings and representation dimension. Adv. Math. 185 (2004), no. 1, 159--177.

[FGR] R. M. Fossum, P. Griffith, I. Reiten: Trivial extensions of abelian categories. Homological algebra of trivial extensions of abelian categories with applications to ring theory. Lecture Notes in Mathematics, Vol. 456. Springer-Verlag, Berlin-New York, 1975.

[FZ1] S. Fomin, A. Zelevinsky: Cluster algebras. I. Foundations. J. Amer. Math. Soc. 15 (2002), no. 2, 497--529.

[FZ2] S. Fomin, A. Zelevinsky: Cluster algebras. II. Finite type classification. Invent. Math. 154 (2003), no. 1, 63--121.

[GLS1] C. Geiss, B. Leclerc, J. Schr\"oer: Semicanonical bases and preprojective algebras. Ann. Sci. Ecole Norm. Sup. (4) 38 (2005), no. 2, 193--253.

[GLS2] C. Geiss, B. Leclerc, J. Schr\"oer: Rigid modules over preprojective algebras, to appear in Invent. Math.

[GS] C. Geiss, J. Schr\"oer: Extension-orthogonal components of preprojective varieties. Trans. Amer. Math. Soc. 357 (2005), no. 5, 1953--1962.

[GN1] S. Goto, K. Nishida: Minimal injective resolutions of Cohen-Macaulay isolated singularities. Arch. Math. (Basel) 73 (1999), no. 4, 249--255. 


[GN2] S. Goto, K. Nishida: Towards a theory of Bass numbers with application to Gorenstein algebras. Colloq. Math. 91 (2002), no. 2, 191--253. 

[G] X. Guo: Representation dimension: an invariant under stable equivalence. Trans. Amer. Math. Soc. 357 (2005), no. 8, 3255--3263.

[Ha] D. Happel: Triangulated categories in the representation theory of finite-dimensional algebras. London Mathematical Society Lecture Note Series, 119. Cambridge University Press, Cambridge, 1988.

[H] M. Hoshino: On dominant dimension of Noetherian rings. Osaka J. Math. 26 (1989), no. 2, 275--280.

[Ho] T. Holm: Representation dimension of some tame blocks of finite groups. Proceedings of the First Sino-German Workshop on Representation Theory and Finite Simple Groups (Beijing, 2002). Algebra Colloq. 10 (2003), no. 3, 275--284.


[I1] O. Iyama: Finiteness of Representation dimension, Proc. Amer. Math. Soc. 131 (2003), no. 4, 1011-1014.

[I2] O. Iyama: Symmetry and duality on $n$-Gorenstein rings. J. Algebra 269 (2003), no. 2, 528--535.

[I3] O. Iyama: Representation dimension and Solomon zeta function. Representations of finite dimensional algebras and related topics in Lie theory and geometry, 45--64, Fields Inst. Commun., 40, Amer. Math. Soc., Providence, RI, 2004.

[I4] O. Iyama: $\tau$-categories III: Auslander orders and Auslander-Reiten quivers, Algebr. Represent. Theory 8 (2005), no. 5, 601--619.

[I5] O. Iyama: The relationship between homological properties and representation theoretic realization of artin algebras, Trans. Amer. Math. Soc. 357 (2005), no. 2, 709--734.

[I6] O. Iyama: Rejective subcategories of artin algebras and orders, preprint arXiv:math.RT/0311281.

[I7] O. Iyama: Higher dimensional Auslander-Reiten theory on maximal orthogonal subcategories, to appear in Adv. Math.

[IT] K. Igusa, G. Todorov: On the finitistic global dimension conjecture for Artin algebras. Representations of algebras and related topics, 201--204, Fields Inst. Commun., 45, Amer. Math. Soc., Providence, RI, 2005.

[KK] H. Krause, D. Kussin: Rouquier's theorem on representation dimension, preprint arXiv:math.RT/0505055.

[L] G. Leuschke: Endomorphism rings of finite global dimension, to appear in Canadian J. Math.

[Ma] H. Matsumura: Commutative ring theory. Cambridge Studies in Advanced Mathematics, 8. Cambridge University Press, Cambridge, 1989.

[Mc] J. McKay: Graphs, singularities, and finite groups. The Santa Cruz Conference on Finite Groups (Univ. California, Santa Cruz, Calif., 1979), pp. 183--186, Proc. Sympos. Pure Math., 37, Amer. Math. Soc., Providence, R.I., 1980. 

[Mi] J. Miyachi: Injective resolutions of Noetherian rings and cogenerators. Proc. Amer. Math. Soc. 128 (2000), no. 8, 2233--2242.

[M] Y. Miyashita: Tilting modules of finite projective dimension. Math. Z. 193 (1986), no. 1, 113--146.

[N] K. Nishida: Cohen-Macaulay isolated singularities with a dualizing module, Algebr. Represent. Theory 9 (2006), no. 1, 13--31.

[P] J. A. de la Pena: On the dimension of the module-varieties of tame and wild algebras. Comm. Algebra 19 (1991), no. 6, 1795--1807.

[R] M. Ramras: Maximal orders over regular local rings of dimension two. Trans. Amer. Math. Soc. 142 (1969), 457--479.

[RV] I. Reiten, M. Van den Bergh: Two-dimensional tame and maximal orders of finite representation type. Mem. Amer. Math. Soc. 80 (1989).

[Ri] C. Riedtmann: Algebren, Darstellungskocher, Uberlagerungen und zuruck. Comment. Math. Helv. 55 (1980), no. 2, 199--224.

[Ro1] R. Rouquier: Representation dimension of exterior algebras, to appear in Invent. Math.

[Ro2] R. Rouquier: Dimensions of triangulated categories, preprint arXiv:math.CT/0310134.

[T] H. Tachikawa: Quasi-Frobenius rings and generalizations. Lecture Notes in Mathematics, Vol. 351. Springer-Verlag, Berlin-New York, 1973.

[V1] M. Van den Bergh: Three-dimensional flops and noncommutative rings. Duke Math. J. 122 (2004), no. 3, 423--455.

[V2] M. Van den Bergh: Non-commutative crepant resolutions. The legacy of Niels Henrik Abel, 749--770, Springer, Berlin, 2004.

[X] C. Xi: Representation dimension and quasi-hereditary algebras. Adv. Math. 168 (2002), no. 2, 193--212.

[Y] Y. Yoshino: Cohen-Macaulay modules over Cohen-Macaulay rings. London Mathematical Society Lecture Note Series, 146. Cambridge University Press, Cambridge, 1990. 

\vskip.5em{\footnotesize
{\sc Department of Mathematics, University of Hyogo, Himeji, 671-2201, Japan}

{\it iyama@sci.u-hyogo.ac.jp}

\vskip.5em
Current address: 

{\sc Graduate School of Mathematics, Nagoya University,

Chikusa-ku, Nagoya, 464-8602, Japan}

{\it iyama@math.nagoya-u.ac.jp}}
\end{document}